\documentclass[a4paper,11pt]{article}
\usepackage[utf8]{inputenc}
\usepackage[T1]{fontenc}
\usepackage{amssymb,amsmath,amsthm}
\usepackage{xspace}
\usepackage{graphicx}
\usepackage[english]{babel}
\usepackage{color}
\usepackage[left=3.5cm, right=3.5cm]{geometry}
\usepackage{mathtools}
\usepackage{color}
\usepackage[citecolor = blue]{hyperref}

\newcommand{\eps}{\varepsilon}
\newcommand{\reals}{\mathbb{R}}

\newcommand{\nats}{\mathbb{N}}

\newcommand{\esp}{\mathbb{E}}
\newcommand{\prob}{\mathbb{P}}
\newcommand{\srf}{\Sigma}

\DeclarePairedDelimiter\floor{\lfloor}{\rfloor}
\newcommand\eqlab{\addtocounter{equation}{1}\tag{\theequation}}

\newenvironment{prf}[1][]
{\vskip 2mm  {\it \bf Proof#1. }}{$\Box$ \vskip 2mm}
\newtheorem{thm}{Theorem}
\newtheorem{lm}[thm]{Lemma}
\newtheorem{prop}[thm]{Proposition}

\title{Hole probability for nodal sets of the cut-off Gaussian Free Field}
\author{Alejandro Rivera}
\date{\today}

\begin{document}

\definecolor{shadecolor}{gray}{0.95}
\maketitle

\begin{abstract}
Let $(\Sigma,g)$ be a closed connected surface equipped with a riemannian metric. Let $(\lambda_n)_{n\in\nats}$ and $(\psi_n)_{n\in\nats}$ be the increasing sequence of eigenvalues and the sequence of corresponding $L^2$-normalized eigenfunctions of the laplacian on $\Sigma$. For each $L>0$, we consider $\phi_L=\sum_{0<\lambda_n\leq L}\frac{\xi_n}{\sqrt{\lambda_n}}\psi_n$ where the $\xi_n$ are i.i.d centered gaussians with variance $1$. As $L\rightarrow\infty$, $\phi_L$ converges a.s. to the Gaussian Free Field on $\Sigma$ in the sense of distributions. We first compute the asymptotic behavior of the covariance function for this family of fields as $L\rightarrow\infty$. We then use this result to obtain the asymptotics of the probability that $\phi_L$ is positive on a given open proper subset with smooth boundary. In doing so, we also prove the concentration of the supremum of $\phi_L$ around $\frac{1}{\sqrt{2\pi}}\ln L$.
\end{abstract}

\selectlanguage{english}
\tableofcontents
\clearpage

\section{Introduction}
\subsection{Setting and main results}
In recent years, there have been many developments in the study of random linear combinations of eigenfunctions of the laplacian on a closed manifold. In this paper we consider a different model with strong ties to statistical mechanics, mentioned both in \cite{schr} (Problem 2.4) and \cite{zeld} (equation $(97)$). Let $(\Sigma,g)$ be a smooth compact connected surface equipped with a riemannian metric. Let $\Delta=d^*d$ be the Laplace operator on $\Sigma$ associated to $g$ and $|dV_g|$ the volume density defined by $g$. Note that with our convention, if $\Sigma$ is the flat torus with coordinates $(x_1,x_2)$, $\Delta = -\partial_{x_1}^2-\partial_{x_2}^2$. For each $L>0$, let $U_L$ be the real vector space spanned by the eigenfunctions of $\Delta$ whose eigenvalues are positive and smaller than $L$. Then, $U_L$ is finite dimensional and
\[ (u,v)\mapsto \int_\Sigma g(\nabla u,\nabla v) |dV_g| \]
defines a scalar product on each $U_L$ which induces a gaussian probability distribution on $U_L$. For each $L>0$, let $\phi_L$ be random variable chosen with this distribution. Then, each instance of $\phi_L$ is a smooth function on $\Sigma$. In particular, for each $L>0$, $(\phi_L(x))_{x\in\Sigma}$ defines a gaussian field on $\Sigma$. We will see that the field $\phi_L$  converges $L\rightarrow\infty$ almost surely in the sense of distributions to the Gaussian Free Field on $\Sigma$, a central object in contemporary statistical mechanics. Following \cite{zeld}, we choose to call $\phi_L$ \textbf{the cut-off Gaussian Free Field} on $\Sigma$ or CGFF for short. While  the definition of the CGFF is formally similar to that of the usual cut-off eigenfunction model (see for instance \cite{nico_crit}), it is actually quite different. Indeed, while the cut-off model exhibits a local scale with polynomial correlations, the CGFF has global logarithmic correlations. We will prove that it is actually much closer to the discrete Gaussian Free Field. For this purpose we will combine methods from statistical mechanics and random geometry, thus creating a new interface between the two subjects.\\
\begin{center}
\includegraphics[width=0.49\textwidth]{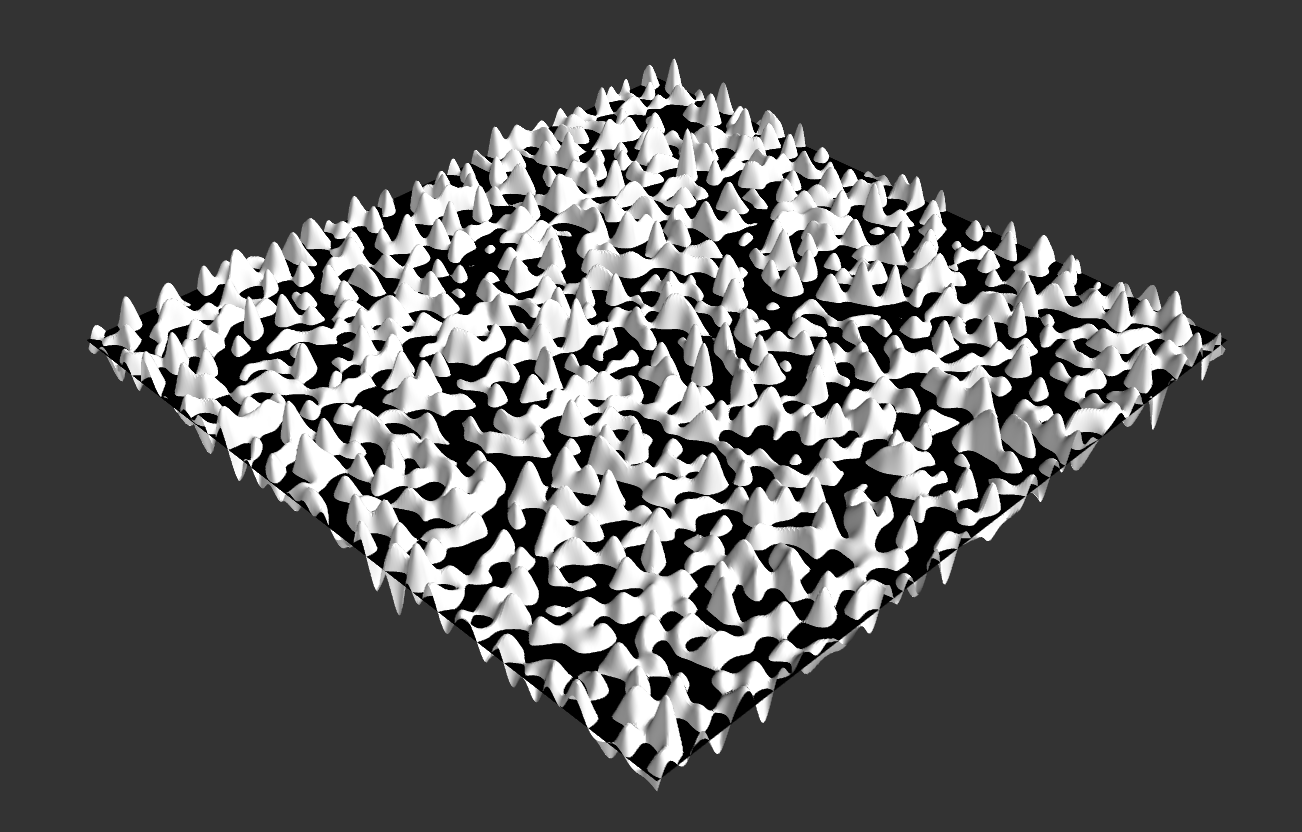}
\includegraphics[width=0.49\textwidth]{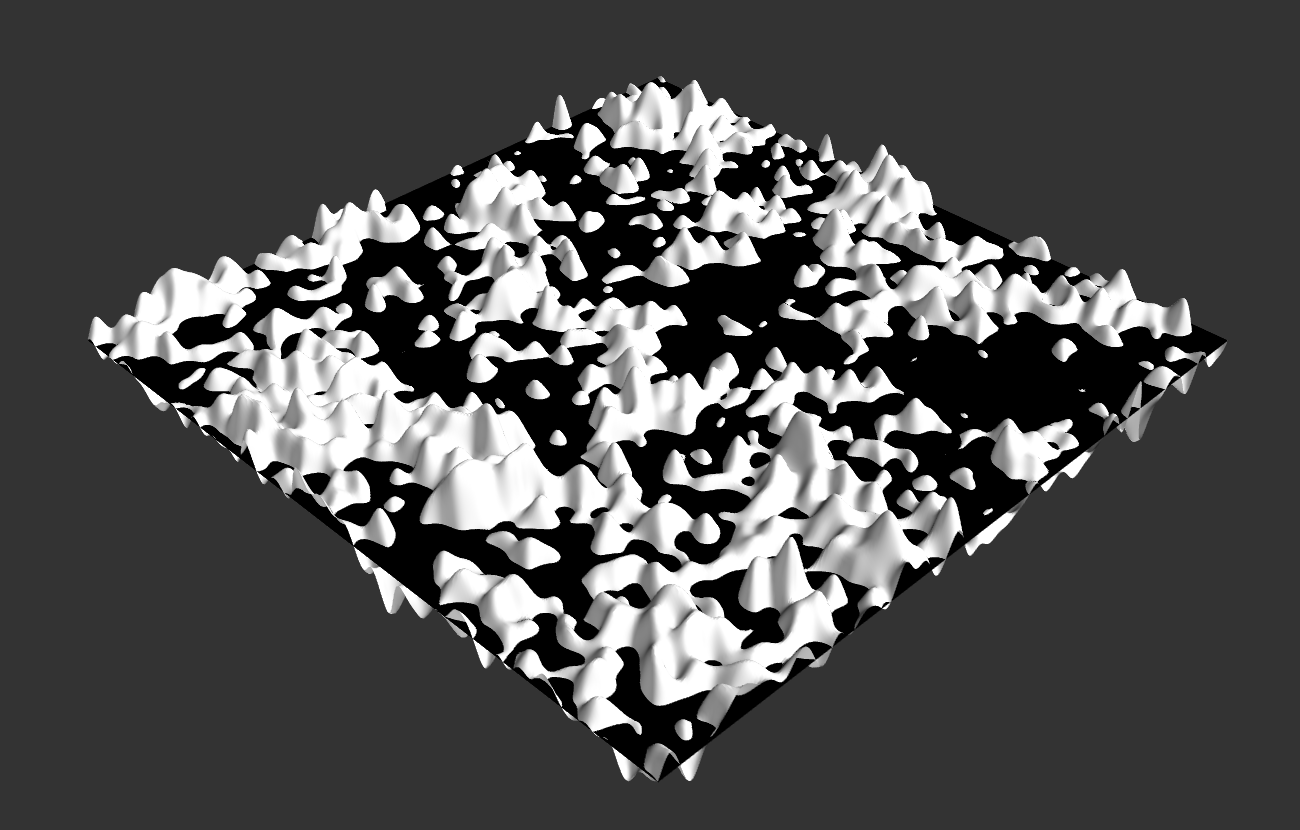}\\
\vspace{5mm}
\textit{An instance of the cut-off field (resp. the CGFF) on the left (resp. right) on the flat torus with $L=1000$. The field is colored in white and the black surface is a horizontal square at height zero.}
\vspace{5mm}
\end{center}
Let $D\subset\Sigma$ be a non-empty proper open subset of $\Sigma$ with smooth boundary. We ask what the probability is that the the field stays positive on $D$. The asymptotic behavior of this probability as $L\rightarrow\infty$ will be expressed in terms of \textbf{the capacity} of $D$, which we define as the infimum of the quantity $\frac{1}{2}\|\nabla h\|_2^2$ taken over all $h\in C^\infty(\Sigma)$ with zero mean, such that $\forall x\in D$, $h(x)\geq 1$ and which we denote by $cap_\Sigma(D)$ or $cap(D)$ when there is no ambiguity. More precisely, we prove the following theorem.
\begin{thm}\label{mainthm}
Let $(\Sigma,g)$ be a smooth compact connected surface without boundary and let $(\phi_L)_{L>0}$ be the CGFF on $(\Sigma,g)$. Let $D$ be a non-empty proper open subset of $\Sigma$ with smooth boundary. Then,
\[ \lim_{L\rightarrow\infty}\frac{\ln\big(\prob(\forall x\in D, \phi_L(x)>0)\big)}{\ln^2\big(\sqrt{L}\big)}=-\frac{2}{\pi}cap(D). \]
\end{thm}
Along the way, we will also prove that $cap(D)>0$. The assumption that $\Sigma$ has no boundary follows the tradition of the study of random sums of eigenfunctions. However, \textbf{all of the results presented in this paper stay valid in the case where $\Sigma$ has a boundary} with minor modifications and provided we change the definitions of the CGFF and capacity accordingly (see section \ref{sect_bd} for the corresponding statements). This is especially significant for two reasons. First, from Riemann's mapping theorem, two non-empty simply connected proper subsets of $\mathbb{C}$ are conformally equivalent. Second, in this setting $cap(D)$ will be \textbf{conformally invariant} (see section \ref{sect_bd} for more details). To the best of our knowledge, this is the first time a non-trivial conformal invariant emerges from the asymptotics of random sums of eigenfunctions. The event of staying positive on a given set has been studied before in the case of sections of complex line bundles on complex manifolds (see \cite{szz_hole}) and in the case of the discrete Gaussian Free Field (or DGFF) on a box in the square lattice (see \cite{bdg_entrop}). In \cite{szz_hole}, Shiffman, Zrebiec and Zelditch actually prove much stronger results relying on large deviation estimates that work because the field they consider has exponential decay in correlations. As will be apparent in the statement of Theorem \ref{kerthm}, this is not the case in our model, which, like \cite{bdg_entrop}, has logarithmic correlations. In this article, Bolthausen, Deuschel and Giacomin prove the following result
\begin{thm}\label{bdgthm}
Let $V_N=\{1,\dots,N\}^2$ be a square box in $\mathbb{Z}^2$. Let $\phi_N$ be the discrete Gaussian Free Field on $V_N$ with Dirichlet boundary conditions. Let $D\subset [0,1]^2$ be an open subset with smooth boundary and at positive distance of $\partial [0,1]^2$. Let $D_N$ be the set of points $y\in V_N$ such that $\frac{1}{N}y\in D$. Then,
\[ \lim_{N\rightarrow\infty}\frac{1}{(\ln N)^2}\ln\Big(\prob(\forall x\in D_N, \phi_N(x)\geq 0)\Big)=-\frac{8}{\pi}cap_V(D).\]
\end{thm}
Here, $cap_V(D)$ is the infimum of $\frac{1}{2}\|\nabla h\|_2^2$ over all the $h\in C^\infty(V)$ with compact support in $\mathring{V}$ such that $h\geq 1$ in $D$. In Theorem \ref{mainthm}, the connectedness assumption simplifies the proofs and is not very restrictive since the CGFF is independent between different components. The assumption that $D$ be an open set with smooth  boundary allows us to use classical results concerning the potential $cap(D)$ and is already present in the discrete setting. Finally, since the field we consider has zero mean, it cannot stay positive on $D=\Sigma$ so we assume $D\neq\Sigma$. We keep the square root inside the logarithm in the statement because it emerges from the proof as a more natural scale in this problem. Our approach follows the structure of \cite{bdg_entrop}. However, we consider fields in a continuous setting, and more importantly, unlike the DGFF, the CGFF does not seem to have a Markov property. The relation with \cite{bdg_entrop} as well as the strategies employed to deal with these issues will be explained below. Finally, we need to estimate the covariance function of the field. This step is central in our strategy since it is only through this object that we can manipulate the field. We prove the following theorem, which is significant in its own right.
\begin{thm}\label{kerthm}
Let $(\Sigma,g)$ be a compact riemannian surface without boundary and let $(\phi_L)_{L>0}$ be the CGFF on $(\Sigma,g)$. 
For each $L> 0$ and $p,q\in\Sigma$, let
\[ G_L(p,q)=\esp[\phi_L(p)\phi_L(q)]. \]
Then, there exists $\eps>0$ such that for each $p,q\in \Sigma$ satisfying $d_g(p,q)\leq \eps$ and for each $L>0$,
\[ G_L(p,q)=\frac{1}{2\pi}\Big(\ln\big(\sqrt{L}\big)-\ln_+\big(\sqrt{L}d_g(p,q)\big)\Big) + \rho_L(p,q) \]
where  $\ln_+(a)=\max(\ln(a),0)$, $d_g$ is the riemannian distance and $\rho_L(p,q)$ is bounded uniformly with respect to $p$, $q$ and $L$.
\end{thm}
The proof of this theorem relies on H\"ormander's estimates for the spectral kernel of an elliptic operator in \cite{ho_sfeo}. The result is reminiscent of the well known analogue for the DGFF (see for instance Lemma 2.2 of \cite{bz_maxgff}). Aside from Theorem \ref{kerthm}, one important step in the proof of Theorem \ref{mainthm} is to control the supremum of the field on a given domain. We prove the following theorem.
\begin{thm}\label{supthm}
Let $(\Sigma,g)$ be a smooth compact riemannian surface without boundary and let $(\phi_L)_{L>0}$ be the CGFF on $(\Sigma,g)$. Let $D$ be a non-empty open subset of $\Sigma$. Then for each $\eta >0$,
\begin{equation*}
\limsup_{L\rightarrow\infty}\frac{\ln\Big(\prob\Big(\sup_\Sigma\phi_L>\Big(\sqrt{\frac{2}{\pi}}+\eta\Big)\ln\big(\sqrt{L}\big)\Big)\Big)}{\ln\big(\sqrt{L}\big)}\leq -2\sqrt{2\pi}\eta+O(\eta^2)
\end{equation*}
and there exists $a>0$ such that for $L$ large enough,
\begin{equation*}
\prob\Big(\sup_D\phi_L\leq\Big(\sqrt{\frac{2}{\pi}}-\eta\Big)\ln\big(\sqrt{L}\big)\Big)\leq \exp\Big(-a\ln^2\big(\sqrt{L}\big)\Big).
\end{equation*}
\end{thm}
The maxima of random fields on smooth manifolds have been studied, for instance for holomorphic sections of line bundles on K\"ahler manifolds in \cite{sz_sup} and for another eigenfunction model in \cite{bule}. However, these fields are not log-correlated so the probabilistic arguments employed are quite different. This theorem is an analogue of Theorem 2 in \cite{bdg_entrop} and the proof relies on it. The supremum of general discrete log-correlated fields has been studied in \cite{drz_maxlog}. In the case of the DGFF, as well as a large class of continuous log-correlated fields, the law of the supremum has been studied with much higher precision, see for instance \cite{bz_maxgff}, \cite{bdz_cvlaw}, \cite{bl_poiss}, \cite{ac_14}, \cite{drz_maxlog} and \cite{mad_log}. To obtain such results for the CGFF, one would need more precise estimates for the covariance function. Note that while \cite{ac_14} and \cite{mad_log} deal with continuous fields similar to the CGFF, the left tail estimate in Theorem \ref{supthm} does not appear in these works.\\

The paper is organised as follows. In the rest of this section, we give an outline of the proof and introduce some basic notation in order to give a more concrete definition of the CGFF. In section \ref{sect_max} we prove Theorem \ref{supthm}. In section \ref{sect_hole} we prove Theorem \ref{mainthm}. Section \ref{sect_ana} is dedicated to the proof of the analytical tools used before, most notably, Theorem \ref{kerthm}. In section \ref{sect_bd}, we cover the case where $\Sigma$ has a boundary. In the appendix \ref{sect_apx}, we recall some classical results regarding the laplacian and the capacity.

\subsection{Comparison with the discrete setting}

Part of this paper is written in the spirit of \cite{bdg_entrop} which studies the hole probability of the discrete Gaussian Free Field. In this section we outline our proof strategy with \cite{bdg_entrop} in mind and explain the new ideas introduced to deal with this model. We will use the notations introduced in Theorem \ref{bdgthm}.\\

\begin{center}
\includegraphics[width=0.45\textwidth]{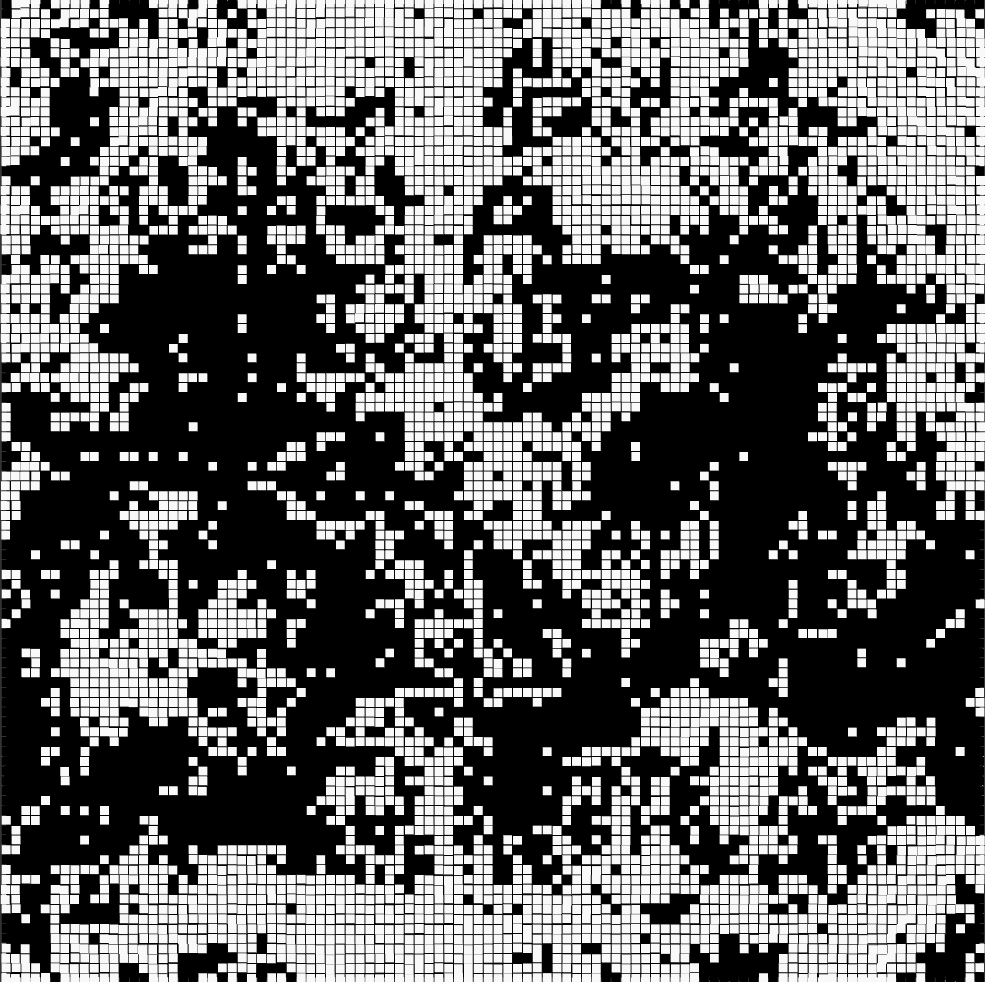}
\hspace{0.05\textwidth}
\includegraphics[width=0.45\textwidth]{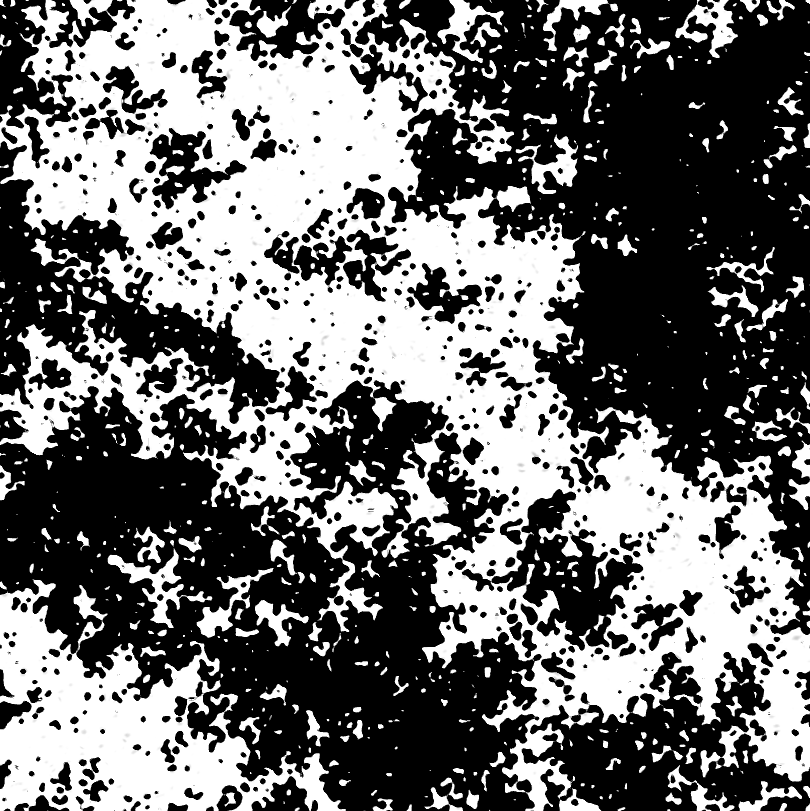}\\
\vspace{5mm}
\textit{On the left, an instance of the DGFF on the square of side $N=100$ with periodic boundary conditions. On the right, an instance of the CGFF on the flat torus with $\sqrt{L}=100$. In both cases, the square is colored white where the field is positive and black where it is negative.}
\vspace{5mm}
\end{center}

To begin with, the field we consider is defined as a random linear combination of eigenfunctions of the laplacian on a compact surface $\Sigma$ and some work is required to obtain a tractable expression for the covariance function. This is Theorem \ref{kerthm} and is proved in section \ref{sect_ana}. In particular we use H\"ormander's estimates of \textbf{the spectral function of the laplacian} from \cite{ho_sfeo}. The analog in the discrete case is taken for granted in \cite{bdg_entrop}. Since this part uses different techniques than the rest of the paper, we present it at the end. We first relate the kernel $G_L$ to the Schwartz kernel of the orthogonal projector onto $U_L$ in order to apply the aforementioned result by H\"ormander. This yields an expression of $G_L$ as an integral involving a kind of generalized Bessel function, which we control thanks to \textbf{the stationary phase method}. It emerges from Theorem \ref{kerthm} that the CGFF $(\phi_L)$ is a log-correlated field that varies at a scale of $L^{-1/2}$. Consequently, in the analogy with \cite{bdg_entrop}, $\Sigma$ will play the role of the square box of size $N$ where $N\simeq L^{1/2}$. This is why we choose to write $\ln\big(\sqrt{L}\big)$ instead of $\frac{1}{2}\ln L$ in our main results above.\\

The proof of Theorem \ref{mainthm} goes as follows. The first step is to estimate the supremum of the field. This is the object of section \ref{sect_max} and the result is Theorem \ref{supthm}. In \cite{bdg_entrop}, the bound on the right tail comes from a simple union bound. In our case, since the space is continuous, the field could fluctuate at scales smaller than $L^{-1/2}$. To control these fluctuations, we use \textbf{a Sobolev inequality} on smalls disks of this scale (see Lemma \ref{sob}). Note that this requires control of the successive derivatives of the field. The result then follows from a union bound applied to a covering of $\Sigma$ by such disks (see Proposition \ref{supup}). To control the left tail of the maximum, Bolthausen, Deuschel and Giacomin use \textbf{the Markov property} of the field to construct a tree-like structure and, inspired by branching random walks, use large deviation results to conclude. At this point we cannot follow the original proof because the CGFF does not seem to have a Markov property. Instead, we restrict the CGFF to a discrete box and use a method from the much more recent \cite{drz_maxlog} in order to apply a gaussian comparison inequality between the restricted field and the DGFF. Thus we recover the bound from the original one (see Proposition \ref{suplow}).\\

Once we have Theorem \ref{supthm}, we can start studying the probability that the CGFF stays positive on a given domain $D$. This is section \ref{sect_hole}. For the lower bound, \cite{bdg_entrop} uses an entropy inequality and the capacity appears by discrete approximation. In our case (see Proposition \ref{poslb}) it seemed more natural to apply \textbf{the barrier method}, already used in \cite{naso_nod} and \cite{gawe_comp}. The idea is to decompose the field into a random multiple of a function $h$ that is greater than one on $D$ and an independent fluctuation. Then, we use the bound on the right tail from Theorem \ref{supthm} to control the supremum of the fluctuation. We then vary $h$ to minimize the cost of this procedure and end up with the capacity of $D$. The lower bound (see Proposition \ref{posub}) is more subtle. Indeed, Bolthausen, Deuschel and Giacomin use the Markov property once again to decompose the DGFF $(\phi_N)$ into two independent gaussian fields. One is ``tamer'' while the other is ``wilder''. We call this \textbf{the two-scale decomposition} of the CGFF. Roughly speaking, to stop the wilder field from making $\phi_N$ negative, the tamer field will need to be close to the expected maximum on a large enough portion of $D$. This will come at a cost that will be related to the capacity. In our paper the Markov property is once again absent. To construct this decomposition, we split the eigenvalue interval in two and obtain a decomposition of $\phi_L$ as an independent sum. The tameness of the tamer field will be immediate. On the other hand the wilder field will require a finer analysis. We will first prove an approximate version of the independence afforded by the Markov property of the DGFF in this decomposition. Then, we will apply once again \textbf{the gaussian comparison method} from \cite{drz_maxlog} to conclude.\\
For the convenience of the reader, we summarize the above discussion in the following table.
\begin{center}
\begin{tabular}{| p{4cm} | p{4cm} | p{4cm} |}
\hline
\textbf{Proof step} & \textbf{Discrete case} & \textbf{Continuous case} \\ \hline
Covariance function &  $\emptyset$ & spectral asymptotics \\ \hline
Right tail for the supremum & union bound & Sobolev inequality + union bound \\ \hline
Left tail of the maximum & Markov property + large deviations & gaussian comparison method \\ \hline
Lower bound for hole probability & entropy inequality & barrier method \\ \hline
Upper bound for hole probability & Markov property + two-scale decomposition & two-scale decomposition + decorrelation estimates + gaussian comparison method\\
\hline
\end{tabular}
\end{center}

\subsection{The cut-off Gaussian Free Field}\label{sect_setting}

From now on, $(\Sigma,g)$ will be a smooth, compact, connected surface without boundary, equipped with a riemannian metric. Let $|dV_g|$ be the density and $\Delta=d^*d$ the Laplace operator defined by $g$. We will denote by $L^2(\Sigma)$ and $H^m(\Sigma)$ for any integer $m\geq 1$ respectively the space of square integrable functions over $\Sigma$ with respect to the measure $|dV_g|$ and the $L^2$ Sobolev space of order $m$ with respect to this same measure (see for instance Definition B.1.1 of \cite{ho_apdo}). For any of these spaces, say $E(\Sigma)$, we will denote by $E_0(\Sigma)$ -- or $E_0$ when no ambiguity is possible -- the subspace of $E(D)$ consisting of functions of zero mean on $\Sigma$. We will denote by $\langle\ ,\ \rangle_2$ the $L^2$ scalar product on $\Sigma$. We will use the same notation in the following case. If $X,Y$ are two vector-fields on $\Sigma$, $\langle X,Y\rangle_2 = \int_\Sigma g_p(X_p,Y_p)|dV_g|(p)$. By the Poincaré-Wirtinger inequality (see Theorem 1 of section 5.8.1 of \cite{evans}) the bilinear form
\[ \langle u,v\rangle_\nabla :=\langle \nabla u,\nabla v\rangle_2. \]
defines a scalar product equivalent to the standard one on $H^1_0(\Sigma)$ called the Dirichlet inner product. It is well known (see for instance Theorem 4.43 of \cite{gahulaf}) that there exist $(\psi_n)_{n\in\nats}\in C^\infty(\Sigma)^\nats$ and $(\lambda_n)_{n\in\nats}\in\reals^\nats$ such that $0=\lambda_0<\lambda_1\leq\lambda_2\dots,\lambda_n\xrightarrow[n\to\infty]{}+\infty$, such that $(\psi_n)_{n\in\nats}$ is a Hilbert basis for $L^2(\Sigma)$ and such that for each $n\in\nats$, $\Delta \psi_n=\lambda_n\psi_n$. In addition, $\psi_0$ is constant and, consequently, $\forall n\geq 1$, $\psi_n\in C^\infty_0(\Sigma)$. Stokes' theorem shows that $\Big(\frac{1}{\sqrt{\lambda_n}}\psi_n\Big)_{n\geq 1}$ is a Hilbert basis of $(H^1_0,\langle\ ,\ \rangle_\nabla)$. For each $L>0$, let $(U_L,\langle\ ,\ \rangle_\nabla)$ be the subspace of $(H^1_0,\langle\ ,\ \rangle_\nabla)$ spanned by the functions $\frac{1}{\sqrt{\lambda_n}}\psi_n$ such that $0<\lambda_n\leq L$. Let $(\xi_n)_{n\in\nats}$ be a sequence of i.i.d real centered gaussian random variables with variance $1$. Then, for each $L>0$ we define the \textbf{cut-off Gaussian Free Field}, or CGFF, as
\begin{equation}\label{sgffeq1}
\phi_L=\sum_{0<\lambda_n\leq L}\frac{\xi_n}{\sqrt{\lambda_n}}\psi_n.
\end{equation}
Hence, for each $L>0$ and $p,q\in\Sigma$,
\begin{equation}\label{kereq1}
G_L(p,q):=\esp[\phi_L(p)\phi_L(q)]=\sum_{0<\lambda_n\leq L}\frac{1}{\lambda_n}\psi_n(p)\psi_n(q).
\end{equation}
Note that defining the CGFF by equation \eqref{sgffeq1} amounts to saying that it is a random function in $U_L$ with probability density proportional to $e^{-\frac{1}{2}\|\phi \|_\nabla^2}d\phi$ where $d\phi$ is the Lebesgue measure on $(U_L,\langle\ ,\ \rangle_\nabla)$. These definitions imply in particular that $\phi_L$ converges almost surely to the Gaussian Free Field in the sense of distributions as $L\rightarrow\infty$ (see for instance section 2.4 of \cite{shef_gff}).

\vspace{2cm}

\textbf{Acknowledgements:} I would like to express my gratitude towards my advisor Damien Gayet for supporting me throughout the course of this project, as well as for his many helpful comments regarding the exposition of these results. I am also grateful to my second advisor Christophe Garban for his encouragements and for sharing his intuitions on the probabilistic aspects of this paper. Finally, I would like to thank Vincent Beffara as well as Luis Alberto Rivera for helping me with the numerical simulations.

\section{The maximum of the CGFF}\label{sect_max}

The aim of this section is to prove Theorem \ref{supthm}. The proof is split in two parts, one for the right tail of the maximum and one for the left. More precisely, Theorem \ref{supthm} follows immediately from Proposition \ref{supup} and Proposition \ref{suplow} below.

\subsection{Binding the right tail the maximum}

The aim of this section is to prove the following proposition.
\begin{prop}\label{supup}
Let $(\Sigma,g)$ be a smooth compact surface and let $(\phi_L)_L$ be the CGFF on $(\Sigma,g)$. Then, for each $\eta>0$, 
\begin{equation*}
\limsup_{L\rightarrow\infty}\frac{\ln\Big(\prob\Big(\sup_\Sigma\phi_L>\Big(\sqrt{\frac{2}{\pi}}+\eta\Big)\ln\big(\sqrt{L}\big)\Big)\Big)}{\ln\big(\sqrt{L}\big)}\leq -2\sqrt{2\pi}\eta+O(\eta^2).
\end{equation*}
\end{prop}
Let us begin by introducing some notation. For each $p\in\Sigma$, $t>0$ and $L>0$, let $D_L(p,t)$ be the riemannian disk of radius $\frac{t}{\sqrt{L}}$ around $p$. For the proof of Proposition \ref{supup}, we will need the following two results.
\begin{lm}\label{sob}
Let $(\Sigma,g)$ be a smooth compact surface and let $(\phi_L)_L$ be the CGFF on $(\Sigma,g)$. Then there is a constant $C>0$ such that for each $p\in\Sigma$ and for $L>0$ large enough,
\begin{align*}
\esp\Big[\sup_{D_L(p,1)}|\phi_L|\Big]&\leq C\sqrt{\ln(L)}\\
\esp\Big[\sup_{D_L(p,1)}|\nabla\phi_L|\Big]&\leq C\sqrt{L}.
\end{align*}
\end{lm}
This lemma is to be compared with Proposition $2.1$ of \cite{gawe_comp}. We postpone its proof till the end of the section. The second result is Theorem 2.1.1 of \cite{adta_rfg} specialised to continuous fields.
\begin{prop}[Borell-Tsirelson-Ibragimov-Sudakov inequality]\label{bor}
Let $T$ be a separable topological space and $(\phi_t)_{t\in T}$ be a centered gaussian field over $T$ which is almost surely bounded and continuous. Then, $\esp[\sup_{t\in T}\phi_t]<\infty$ and for all $u>0$,
\begin{equation*}
\prob\Big(\sup_{t\in T}\phi_t-\esp[\sup_{t\in T}\phi_t]>u\Big)\leq \exp\Big(-\frac{u^2}{2\sigma_T^2}\Big)
\end{equation*}
where $\sigma_T^2=\sup_{t\in T}Var(\phi_t)$.
\end{prop}

Let us now prove Proposition \ref{supup} using Lemma \ref{sob} and Proposition \ref{bor}.

\begin{prf}[ of Proposition \ref{supup}]
According to Theorem \ref{kerthm} for $x=y$, there is a constant $C$ such that for each $p\in\Sigma$,
\[\sup_{q\in D_L(p,1)}Var(\phi_L(q))=\sup_{q\in D_L(p,1)}G_L(q,q)\leq \frac{1}{2\pi}\ln\big(\sqrt{L}\big)+C.\]
Let $\eta>0$. We apply Proposition \ref{bor} to $(\phi_L(q))_{q\in D_L(p,1)}$ using Lemma \ref{sob} to deduce that for each $0<\eta'<\eta$, for $L$ large enough and for all $p\in\Sigma$,
\begin{equation*}
\prob\Big(\sup_{D_L(p,1)}\phi_L>\Big(\sqrt{\frac{2}{\pi}}+\eta\Big)\ln\big(\sqrt{L}\big)\Big)\leq \exp\Big(-\Big(2+2\sqrt{2\pi}\eta'+O(\eta^{'2})\Big)\ln\big(\sqrt{L}\big)\Big).
\end{equation*}
Choose some $0<\eta'<\eta$. Since $\Sigma$ is compact, there exists another constant which we also denote by $C$, such that for each $L$, $\Sigma$ is covered by $CL$ disks of radius $\frac{1}{\sqrt{L}}$. Thus,
\begin{align*}
\prob\Big(\sup_\Sigma\phi_L>\Big(\sqrt{\frac{2}{\pi}}+\eta\Big)\ln\big(\sqrt{L}\big)\Big)&\leq CL\exp\Big(-\Big(2+2\sqrt{2\pi}\eta'+O(\eta^{'2})\Big)\ln\big(\sqrt{L}\big)\Big)\\
                                                                  &\leq C\big(\sqrt{L}\big)^{-2\sqrt{2\pi}\eta'+O(\eta^{'2})}.
\end{align*}
Since $-\phi_L$ has the same law as $\phi_L$, we have the analogous result for the minimum. Therefore, for each $0<\eta'<\eta$ there is an $L_0$ such that for each $L\geq L_0$,
\begin{equation*} 
\prob\Big(\sup_\Sigma|\phi_L|>\Big(\sqrt{\frac{2}{\pi}}+\eta\Big)\ln\big(\sqrt{L}\big)\Big)\leq \big(\sqrt{L}\big)^{-2\sqrt{2\pi}\eta'+O(\eta'^2)}.
\end{equation*}
\end{prf}

To prove Lemma \ref{sob}, we use Theorem \ref{kerthm} and the following proposition. The proof of both of these results is presented in the last section.

\begin{prop}\label{kerprop}
Let $Q_1$ and $Q_2$ be differential operators on $\Sigma$ of respective orders $d_1$ and $d_2$ and let $d=d_1+d_2$. Suppose that $d\geq 1$. Then there exists $C>0$ such that for each $p\in\Sigma$ and $L>0$,
\begin{equation*}
\Big|(Q_1\otimes Q_2)G_L(p,p)\Big|\leq C\big(1+L^{d/2}\big).
\end{equation*}
\end{prop}

\begin{prf}[ of Lemma \ref{sob}]
Let $p\in\srf$. We apply the Sobolev inequality from paragraph 5.2.4 of \cite{fe_gmt} with $m=2$ and $N=2$. The inequality implies there exist constants $C,L_0>0$ such that for all $L\geq L_0$,
\begin{align*}
\sup_{q\in D_L(p,1)}|\phi_L(q)|\leq C\Big[ &\Big(\frac{1}{Vol(D_L(p,2))}\int_{D_L(p,2)}\phi_L(q)^2|dV_g|(q)\Big)^{\frac 1 2}\\
                                                                           &+L^{-\frac 1 2}\Big(\frac{1}{Vol(D_L(p,2))}\int_{D_L(p,2)}|\nabla \phi_L(q)|^2|dV_g|(q)\Big)^{\frac 1 2}\\
                                                                           &+L^{-1}\Big(\frac{1}{Vol(D_L(p,2))}\int_{D_L(p,2)}|\nabla^2\phi_L(q)|^2|dV_g|(q)\Big)^{\frac 1 2}\Big].
\end{align*}
Here $\nabla^2$ denotes the hessian defined by the metric $g$. Note that in order to bind the supremum of a function on a two-dimensional space by $L^2$ Sobolev norms, one must use derivatives up to order at least two. By compactness, $C$ and $L_0$ may be chosen independent of $p$. The same inequality holds for expectations and applying Jensen's inequality to the right hand side we obtain :
\begin{align*}
\esp[\sup_{q\in D_L(p,1)}|\phi_L(q)|]\leq C\Big[ &\Big(\frac{1}{Vol(D_L(p,2))}\int_{D_L(p,2)}\esp[\phi_L(q)^2]|dV_g|(q)\Big)^{\frac 1 2}\\
                                                                                 &+L^{-\frac 1 2}\Big(\frac{1}{Vol(D_L(p,2))}\int_{D_L(p,2)}\esp[|\nabla \phi_L(q)|^2|]dV_g|(q)\Big)^{\frac 1 2}\\
                                                                                 &+L^{-1}\Big(\frac{1}{Vol(D_L(p,2))}\int_{D_L(p,2)}\esp[|\nabla^2\phi_L(q)|^2]|dV_g|(q)\Big)^{\frac 1 2}\Big].
\end{align*}
In the above inequality, for any tensor $T$, $|T|$ denotes the norm of $T$ induced by $g$ on the corresponding tensor bundle. Since for any differential operator $P$ over $\srf$, $\esp[P\phi_L(q)^2]=(P\otimes P)G_L(q,q)$, we get
\begin{align*}
\esp[\sup_{q\in D_L(p,1)}|\phi_L(q)|]\leq C\Big[ &\Big(\frac{1}{Vol(D_L(p,2))}\int_{D_L(p,2)}G_L(q,q)|dV_g|(q)\Big)^{\frac 1 2}\\
                                                             &+L^{-\frac 1 2}\Big(\frac{1}{Vol(D_L(p,2))}\int_{D_L(p,2)}|(\nabla\otimes\nabla)G_L(q,q)|\ |dV_g|(q)\Big)^{\frac 1 2}\\
                                                             &+L^{-1}\Big(\frac{1}{Vol(D_L(p,2))}\int_{D_L(p,2)}|(\nabla^2\otimes\nabla^2)G_L(q,q)|\ |dV_g|(q)\Big)^{\frac 1 2}\Big]. \eqlab \label{ref1}
\end{align*}
Now, from Proposition \ref{kerprop}, there is a constant $C>0$ such that for all $p\in\srf$ and $L>0$,
\[ |G_L(p,p)|\leq C\ln(L);\hspace{5mm} |(\nabla\otimes\nabla)G_L(p,p)|\leq C L;\hspace{5mm} |(\nabla^2\otimes\nabla^2)G_L(p,p)|\leq C L^2. \]
Applying these three inequalities to equation \eqref{ref1}, we deduce that there is a constant $C>0$ such that for all $p\in\srf$ and $L>0$,
\begin{equation*}
\esp[\sup_{q\in D_L(p,1)}|\phi_L(q)|]\leq C\sqrt{\ln(L)}.
\end{equation*}
This proves the first statement. The proof carries over to the second statement almost verbatim, using in addition the following estimate from Proposition \ref{kerprop}.
\begin{equation*}
|(\nabla^3\otimes\nabla^3)G_L(p,p)|\leq C L^3.
\end{equation*}
\end{prf}

\subsection{Binding the left tail of the maximum}

In this section we prove the following proposition.
\begin{prop}\label{suplow}
Let $(\Sigma,g)$ be a smooth compact surface and let $(\phi_L)_L$ be the CGFF on $(\Sigma,g)$. Let $D\subset \Sigma$ be a non-empty open subset of $\Sigma$. Then, for each $\eta>0$ there is a constant $a>0$ such that for $L$ large enough,
\begin{equation*}
\prob\Big(\sup_D\phi_L\leq\Big(\sqrt{\frac{2}{\pi}}-\eta\Big)\ln\big(\sqrt{L}\big)\Big)\leq \exp\Big(-a\ln^2\big(\sqrt{L}\big)\Big).
\end{equation*}
\end{prop}

We now introduce some notation. For each $N\in\nats_{\geq 1}$, let $V_N$ be the set of points $(x_1,x_2)\in\mathbb{Z}^2$ such that $0\leq x_1,x_2 \leq N-1$, let $V_N'$ the set of points in $V_N$ at distance at least $N/4$ from the boundary and let $x_N$ be one of the points nearest to its center. For each $x,y\in V_N$, $|x-y|$ will denote the euclidian distance between $x$ and $y$. For each $t>0$, $\ln_+t$ will denote $\max(\ln t,0)$.\\

Proposition \ref{suplow} will follow from the two following results.
\begin{prop}\label{lcmax}
Let $(X_t)_{t>1}$ be a family of random fields such that for all $t>1$, $X_t$ is defined over the box $V_{\floor{t}}$. Suppose there is a constant $C>0$ such that for each $t>1$ and $x,y\in V_{\floor{t}}$,
\begin{equation*}
\big|\esp[X_t(x)X_t(y)]-\ln t +\ln_+|x-y|\big|\leq C.
\end{equation*}
Then for each $\eta>0$ there is a constant $a>0$ depending only on $C$ and $\eta$ such that for $t>1$ large enough,
\begin{equation*}
\prob\Big(\sup_{V_{\floor{t}}}X_t\leq (2-\eta)\ln t\Big)\leq\exp\big(-a(\ln t)^2\big).
\end{equation*}
\end{prop}
\begin{lm}\label{cast0}
Fix $0<\delta<\frac{1}{2\sqrt{2}}$ and $L>0$. Let $\iota : V_{\floor{\sqrt{L}}}\rightarrow\Sigma$ be an injection of the $\mathbb{Z}^2$-box of side-length $\floor{\sqrt{L}}$ into $\Sigma$ such that for any distinct $x,y\in V_{\floor{\sqrt{L}}}$,
\[\frac{\delta}{2}\leq\frac{\sqrt{L}d_g(\iota(x),\iota(y))}{|x-y|}\leq 2\delta. \]
Then, there is a constant $C(\delta)>0$ independent of $L$ and $\iota$ such that for $L$ large enough and for each $x,y\in V_{\floor{\sqrt{L}}}$,
\[\Big|\esp[\phi_L(\iota(x))\phi_L(\iota(y))]-\frac{1}{2\pi}\ln\big(\sqrt{L}\big)+\ln_+|x-y|\Big|\leq C(\delta). \]
\end{lm}
Lemma \ref{cast0} is just the specialisation of Lemma \ref{cast} from section \ref{sect_hole} to the case $\alpha=0$. In the following proof we use Proposition \ref{lcmax} and Lemma \ref{cast0}.
 
\begin{prf}[ of Proposition \ref{suplow}]
Choose some $\delta>0$ and $\iota$ satisfying the properties required to apply Lemma \ref{cast0} and such that the image of $\iota$ is contained in $D$ for $L$ large enough. Then, by Lemma \ref{cast0}, the family $(X_t)_t$ defined by, $\forall t>0$, $X_t:=\sqrt{2\pi}\phi_{t^2}\circ\iota$ satisfies the hypotheses of Proposition \ref{lcmax}. In particular, for each $\eta>0$, there is a constant $a>0$ such that for all $L>0$,
\begin{align*}
\prob(\sup_D\phi_L\leq\Big(\sqrt{\frac{2}{\pi}}-\eta\Big)\ln\big(\sqrt{L}\big)\Big)&\leq\prob\Big(\sup_{V_{\floor{\sqrt{L}}}}X_{\sqrt{L}}\leq(2-\sqrt{2\pi}\eta)\ln\big(\sqrt{L}\big)\Big)\\
                                                                                   &\leq \exp\Big(-a\ln^2\big(\sqrt{L}\big)\Big).
\end{align*}
\end{prf}

To prove Proposition \ref{lcmax}, we use the two following results. The first is a special case of Theorem 2 (b) of \cite{bdg_entrop}.
\begin{thm}\label{bdgmaxlb}
Let $(\phi_N)_N$ be the DGFF on $V_N$ with the standard normalization. For each $\eta>0$ there is a constant $c>0$ such that for each $N\in\nats$ large enough,
\begin{equation*}
\prob\Big(\sup_{V_N'}\phi_N\leq\Big(\sqrt{\frac{8}{\pi}}-\eta\Big)\ln N\Big)\leq \exp(-c(\ln N)^2).
\end{equation*}
\end{thm}
\begin{lm}[Slepian's Lemma, see Theorem 2.2.1 of \cite{adta_rfg}]\label{slep}
Let $T$ be a separable topological space and $(Z(p))_{p\in T}$, $(Y(p))_{p\in T}$ be two continuous centered gaussian fields on $T$ satisfying the following two properties.
\begin{enumerate}
\item For each $p\in T$, $\esp[Z(p)^2]=\esp[Y(p)^2]$.
\item For each $p,q\in T$, $\esp[Z(p)Z(q)]\leq\esp[Y(p)Y(q)]$.
\end{enumerate}
Then, for each $u\in\reals$,
\begin{equation*}
\prob(\sup_{p\in T}Z(p)>u)\geq\prob(\sup_{p\in T}Y(p)>u).
\end{equation*}
\end{lm}
We now deduce Proposition \ref{lcmax} from Proposition \ref{bdgmaxlb} and Lemma \ref{slep}. The following proof is inspired by that of Lemma 2.8 of \cite{drz_maxlog}.
\begin{prf}[ of Proposition \ref{lcmax}]
Choose any $\eta>0$ and some $j\in\nats$ to be fixed later. For each $t>2^j$, let $N=N(t)=\floor{2^{-j}t}$ and let $Z_{N}$ be gaussian field defined on $V_N$ by setting, for each $x\in V_N$, $Z_N(x)=X_t(2^jx)$. Then, for distinct $x,y\in V_N$,
\begin{align*}
\big|\esp[Z_N(x)^2]-\ln N -j\ln 2\big|&\leq C\\
\big|\esp[Z_N(x)Z_N(y)]-\ln N + \ln_+|x-y|\big|&\leq C.
\end{align*}
Let $Y_N$ be the DGFF on $V_N$ multiplied by $\sqrt{\frac{\pi}{2}}$. From Lemma 2.2 of \cite{bz_maxgff} there is a universal constant $C_0>0$ such that for any distinct $x,y\in V_N'$,
\begin{align*}
\big|\esp[Y_N(x)^2]-\ln N\big|&\leq C_0\\
\big|\esp[Y_N(x)Y_N(y)]-\ln N + \ln_+|x-y|\big|&\leq C_0.
\end{align*}
Since $1/4<\ln 2<1$, there is $j_0$ depending only on $C$ (and the universal constant $C_0$) such that for each $j\geq j_0$ and for each $x\in V_N'$,
\begin{equation*}
j/4\leq\esp[Z_N(x)^2]-\esp[Y_N(x)^2]\leq j.
\end{equation*}
Let
\begin{equation*}
a_N(x)=\sqrt{j^{-1}(\esp[Z_N(x)^2]-\esp[Y_N(x)^2])}\in[1/2,1]
\end{equation*}
and choose $\xi$ a centered gaussian random variable with variance $1$ independent from the fields previously introduced. Then, there is $j\geq j_0$ depending only on $C$ such that for each $x,y\in V_N'$ distinct,
\begin{align*}
\esp[Z_N(x)^2]&=\esp[(Y_N(x)+\sqrt{j}\xi a_N(x))^2]\\
\esp[Z_N(x)Z_N(y)]&\leq\esp[(Y_N(x)+\sqrt{j}\xi a_N(x))(Y_N(y)+\sqrt{j}\xi a_N(y))].
\end{align*}
Thus, by Lemma \ref{slep},
\begin{align*}
\prob\Big(\sup_{V_{\floor{t}}}X_t\leq (2-\eta)\ln t\Big)\leq&\prob\Big(\sup_{V_N'}Z_N\leq (2-\eta)\ln t\Big)\\
                                                   \leq&\prob\Big(\sup_{x\in V_N'}\Big[Y_N(x)+\sqrt{j}\xi a_N(x)\Big]\leq (2-\eta)\ln t\Big)\\
                                                   \leq&\prob\Big(\sup_{x\in V_N'}Y_N(x)\leq (2-(\eta/2))\ln t\Big)+\\
                                                       &\prob\Big(\xi\geq \frac{\eta}{2\sqrt{j}\sup_{V_N'}a_N}\ln t\Big).
\end{align*}
For $t$ large enough, $(2-(\eta/2))\ln t\leq (2-(\eta/3))\ln(N)$. From standard tail estimates for gaussian variables applied to $\xi$ and Theorem \ref{bdgmaxlb} applied to $\sqrt{\frac{2}{\pi}}Y_N$, there is a constant $a>0$ such that for $t$ large enough,
\begin{equation*}
\prob\Big(\sup_{V_{\floor{t}}'}X_t\leq (2-\eta)\ln t\Big)\leq\exp\big(-a(\ln t)^2\big).
\end{equation*}
Moreover, $a$ depends only on $C$, $\eta$.
\end{prf}

\section{Hole probabilitiy for the CGFF}\label{sect_hole}

The aim of this section is to prove Theorem \ref{mainthm}. From now on, we fix an open proper subset $D$ of $\Sigma$ with smooth boundary. This assumption implies that there exist functions $f\in C^\infty_0(\Sigma)$ greater or equal to $1$ on $D$. We want to estimate the probability of the event
\begin{equation*}
\Omega_L^+=\Big\{\forall x\in D, \phi_L(x)>0\Big\}.
\end{equation*}
We divide the proof into two parts, one for the lower bound and one for the upper bound. Theorem \ref{mainthm} will thus follow from immediately from Proposition \ref{poslb} and Proposition \ref{posub} below. Note that, if $D$ is non-empty, then Proposition \ref{capalt} implies $cap(D)>0$.

\subsection{The two-scale decomposition}

In this section, we introduce the two-scale decomposition of the CGFF used below. We will use notations from sections \ref{sect_setting} and \ref{sect_max}. For each $\alpha\in ]0,1]$ and each $L>0$, we denote by $\psi_{\alpha,L}$ the field $\psi_{\alpha,L}=\sum_{L^\alpha<\lambda_n\leq L}\frac{\xi_n}{\sqrt{\lambda_n}}\psi_n$. Note that
\begin{equation}\eqlab \label{ref4}
\phi_L = \psi_{\alpha,L} + \phi_{L^\alpha}
\end{equation}
and that $\psi_{\alpha,L}$ and $\phi_{L^\alpha}$ are independent. Moreover, observe that the two point correlation function of $\psi_{\alpha,L}$ is $G_{L^\alpha,L}(p,q) = \sum_{L^\alpha<\lambda_n\leq L}\frac{1}{\lambda_n}\psi_n(p)\psi_n(q)=G_L(p,q)-G_{L^\alpha}(p,q)$. The asymptotics of $G_{L^\alpha,L}$ will follow easily from those of $G_L$. The field $\phi_{L^\alpha}$ will vary at scale $L^{-\alpha/2}$ while $\psi_{\alpha,L}$ will vary at scale $L^{-1/2}$ and will decorrelate at large distances. While the first fact is immediate, the other two require some justification.\\
\begin{center}
\includegraphics[width=0.4\textwidth]{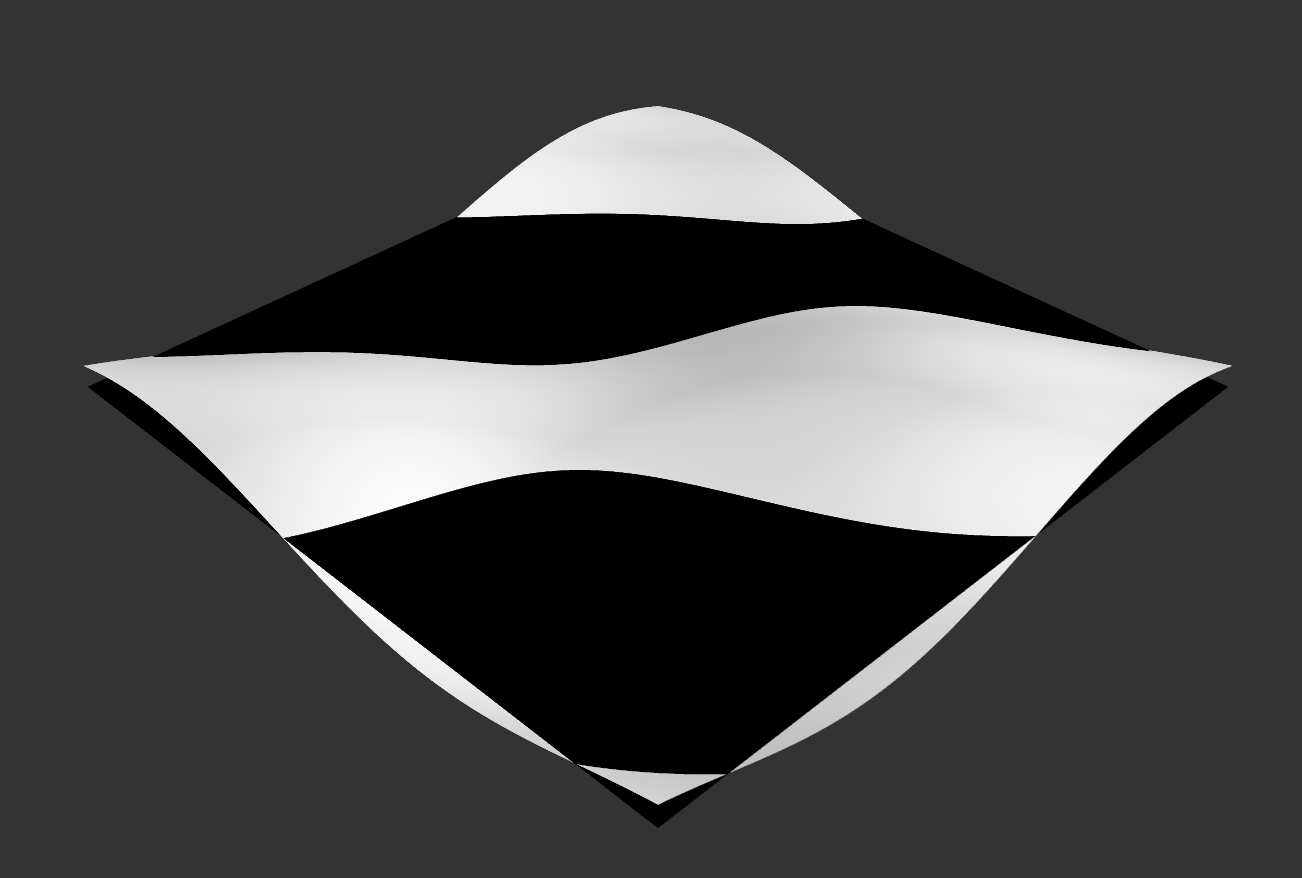}
\includegraphics[width=0.4\textwidth]{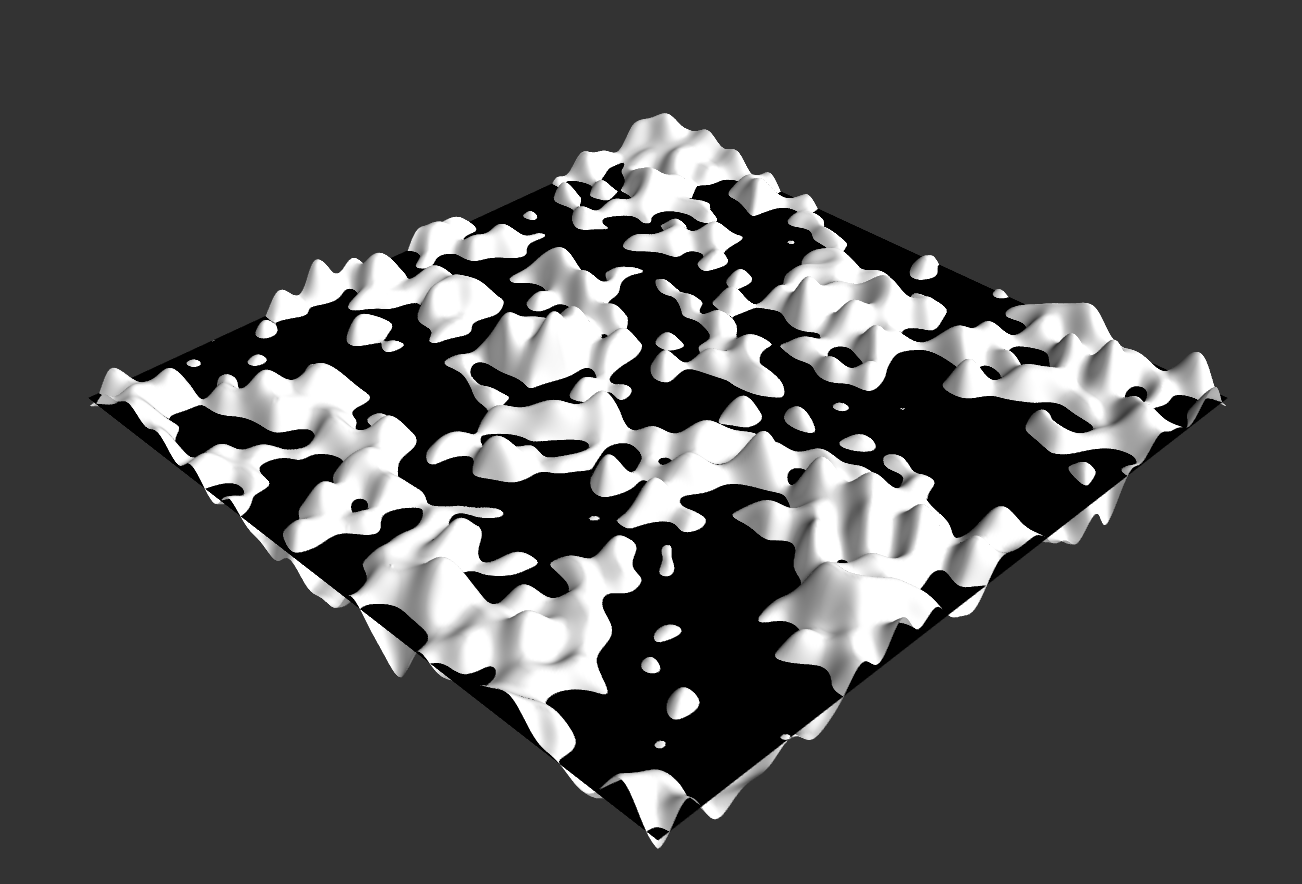}\\
\vspace{5mm}
\textit{An instance of the field $\phi_{L^\alpha}$ on the left and of the field $\psi_{\alpha,L}$ on the right for $\alpha = 0.25$ and $L=400$ on the flat torus. The fields are colored in white and the black surface is a horizontal square at height zero.}\\
\vspace{5mm}
\includegraphics[width=0.5\textwidth]{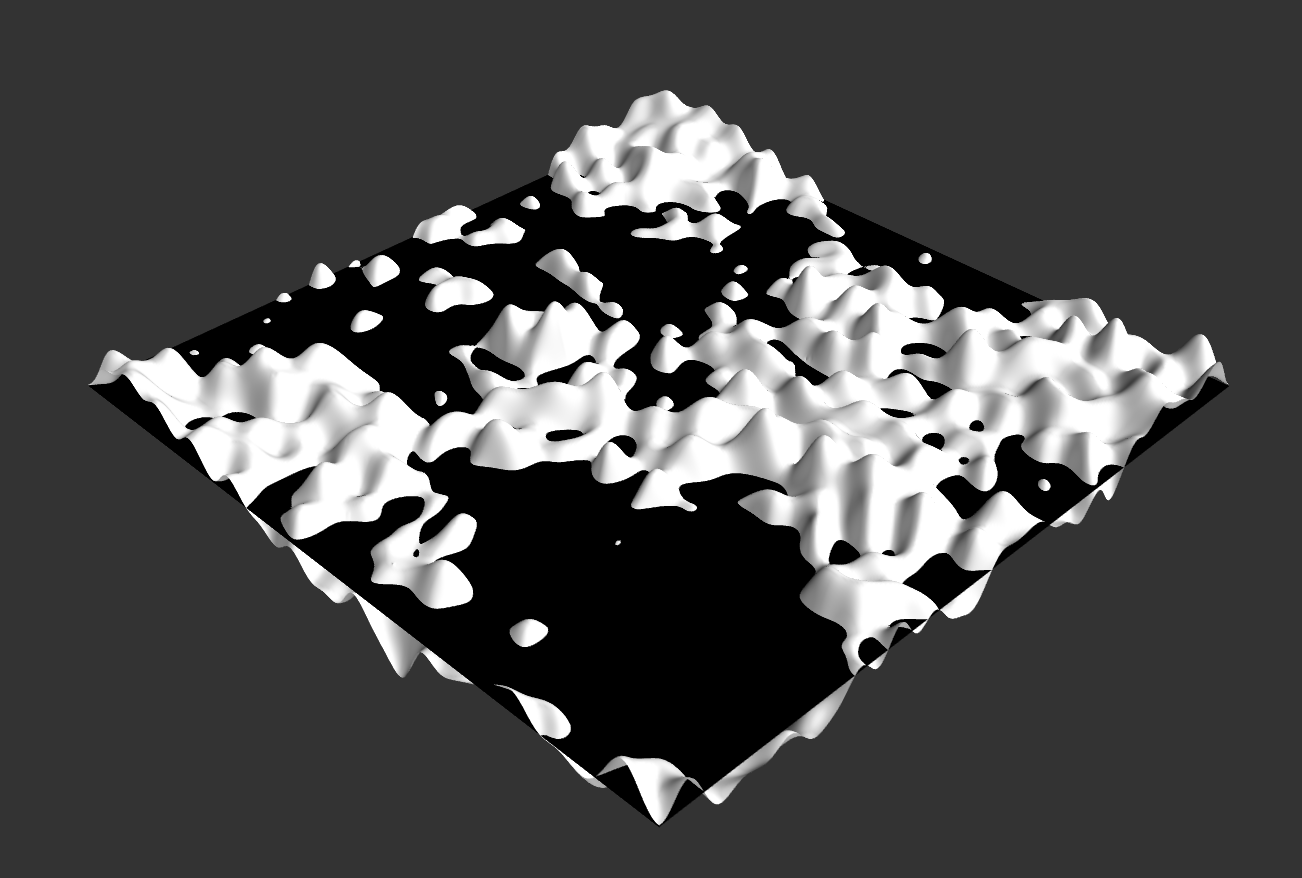}\\
\vspace{5mm}
\textit{The total field $\phi=\phi_{L^\alpha}+\psi_{\alpha,L}$ obtained from the instances above.}
\vspace{5mm}
\end{center}
We will prove the following two results.

\begin{lm}\label{cast}
Fix $0<\delta<\frac{1}{2\sqrt{2}}$, $\alpha\in[0,1[$ and $L>0$. Let $\iota:V_{\floor{L^{(1-\alpha)/2}}}\rightarrow\Sigma$ be an injection of the $\mathbb{Z}^2$-box of side-length $\floor{L^{(1-\alpha)/2}}$ into $\Sigma$ such that for distinct $x,y\in V_{\floor{L^{(1-\alpha)/2}}}$,
\begin{equation*}
\frac{\delta}{2}\leq\frac{\sqrt{L}d_g(\iota(x),\iota(y))}{|x-y|}\leq 2\delta.
\end{equation*}
Then, there is a constant $C(\delta)>0$ independent of $\alpha$, $L$ and $\iota$ such that for $L$ large enough and for each $x,y\in V_{\floor{L^{(1-\alpha)/2}}}$,
\begin{equation*}
\big|\esp[\psi_{\alpha,L}(\iota(x))\psi_{\alpha,L}(\iota(y))]-\frac{1-\alpha}{2\pi}\ln\big(\sqrt{L}\big)+\ln_+|x-y|\big|\leq C(\delta).
\end{equation*}
\end{lm}

This lemma shows that $\psi_{\alpha,L}$ does indeed vary at scale $L^{-1/2}$ and we will use it to prove that its maximum on a box of sidelength $L^{\alpha/2}$ will be close to $\frac{1-\alpha}{2\pi}\ln\big(\sqrt{L}\big)$. Note that Lemma \ref{cast0} is just Lemma \ref{cast} with $\alpha=0$ as announced above. Before proving Lemma \ref{cast}, let us state the second result we will need concerning $\psi_{\alpha,L}$.

\begin{prop}\label{indep}
Choose $0<\alpha<1$. For each $\delta>0$ and $0<\beta<\alpha$,
\[G_{L^\alpha,L}(p,q)\xrightarrow[L\to\infty]{} 0. \]
\end{prop}

This proposition shows that the field decorrelates at large distances. The proof is completely analytical so we leave it for section \ref{sect_ana}.

\begin{prf}[ of Lemma \ref{cast}]
For $L$ large enough, then $V_L:=\iota\big(V_{\floor{L^{(1-\alpha)/2}}}\big)$ has diameter smaller than the $\eps$ in the statement of Theorem \ref{kerthm}. For each $x,y\in V_{\floor{L^{(1-\alpha)/2}}}$,
\begin{align*}
\esp[\psi_{\alpha,L}(\iota(x))\psi_{\alpha,L}(\iota(y))]&=G_{L^\alpha,L}(\iota(x),\iota(y))\\
                                                                     &=G_L(\iota(x),\iota(y))-G_{L^\alpha}(\iota(x),\iota(y)).
\end{align*}
By Theorem \ref{kerthm} applied for $L'=L$ or $L'=L^\alpha$, there is a constant $C>0$ such that for each $L>0$ large enough, for each $x,y\in V_{\floor{L^{(1-\alpha)/2}}}$,
\[ G_L(\iota(x),\iota(y))=\frac{1}{2\pi}\Big(\ln\big(\sqrt{L}\big)-\ln_+\big(\sqrt{L}d_g(\iota(x),\iota(y))\big)\Big)+\rho^1_L(x,y). \]
and
\begin{align*}
G_{L^\alpha}(\iota(x),\iota(y))&=\frac{1}{2\pi}\Big(\ln\big(\sqrt{L^\alpha}\big)-\ln_+\big(L^{\alpha/2}d_g(\iota(x),\iota(y))\big)\Big)+\rho^2_L(x,y)\\
                            &=\frac{1}{2\pi}\Big(\alpha\ln\big(\sqrt{L}\big)-\ln_+\big(L^{\alpha/2}d_g\big(\iota(x),\iota(y))\big)\Big)+\rho^2_L(x,y).
\end{align*}
where $|\rho^j_L(x,y)|\leq C$ for $j=1,2$. Therefore,
\begin{equation*}
\Big|\esp[\psi_L(\iota(x))\psi_L(\iota(y))]-\frac{1-\alpha}{2\pi}\ln\big(\sqrt{L}\big)+\ln_+\big(\sqrt{L}d_g(\iota(x),\iota(y))-\ln_+\big(L^{\alpha/2}d_g\big(\iota(x),\iota(y))\big)\Big|\leq 2C.
\end{equation*}
For each $x,y\in V_{\floor{L^{(1-\alpha)/2}}}$, $|x-y|\leq \sqrt{2}\floor{L^{(1-\alpha)/2}}$ so that, since $\delta<\frac{1}{2\sqrt{2}}$,
\[ L^{\alpha/2}d_g(\iota(x),\iota(y))\leq 2\delta L^{(\alpha-1)/2}|x-y|<1. \]
Therefore, $\ln_+(L^{\alpha/2}|d_g(\iota(x),\iota(y))=0$. Now, for each $x,y\in V_{\floor{L^{(1-\alpha)/2}}}$,
\[ \frac{1}{2}\delta|x-y|\leq \sqrt{L}d_g(\iota(x),\iota(y))\leq 2\delta|x-y| \]
so that 
\[ \ln_+|x-y| -\ln 2 +\ln(\delta) \leq \ln_+\big(\sqrt{L}d_g(\iota(x),\iota(y))\big)\leq \ln_+|x-y| + \ln 2 + \ln(\delta). \]
This concludes the proof of the lemma.
\end{prf}

\subsection{The lower bound}

In this section we use the upper bound in Theorem \ref{supthm} to prove the lower bound in Theorem \ref{mainthm}. In other words, we will prove the following.
\begin{prop}\label{poslb}
Let $(\Sigma,g)$ be a compact smooth surface equipped with a riemannian metric and $(\phi_L)$ be the CGFF on $\Sigma$. Let $D$ be a proper open subset of $\Sigma$ and $\Omega_L^+$ be the event that $\phi_L(x)>0$ for each $x\in D$. Then,
\begin{equation*}
\liminf_{L\rightarrow\infty}\frac{\ln\prob(\Omega_L^+)}{\ln^2\big(\sqrt{L}\big)}\geq -\frac{2}{\pi}cap(D).
\end{equation*}
\end{prop}
Our approach in the following proof is inspired by that of Nazarov and Sodin in section 3 of \cite{naso_nod} and that of Gayet and Welschinger in section 2.2 of \cite{gawe_comp}.
\begin{prf}[ of Proposition \ref{poslb}]
Let us choose $\eps>0$ and a function $h\in C^\infty_0(\Sigma)$ such that for each $x\in D$, $h(x)\geq 1$. Since $u\mapsto ||\nabla u||_2$ is $C^1$ continuous, by Lemma \ref{bar_dens}, for $L$ large than some $L_0$, there is a function $f\in U_L$ such that $\forall x\in \Sigma$, $|f(x)-h(x)|\leq \eps$ and such that $\|\nabla f-\nabla h\|_2\leq \eps$. Now, for $L$ large enough, the random field $\phi_L$ can be decomposed as the independent sum $\xi \frac{f}{\|\nabla f\|_2}+\tilde{\phi}_L$ where $\xi$ is a real centered gaussian random variable with variance $1$ and $\tilde{\phi}_L$ is some gaussian field. Choose $\tilde{\xi}$ another real centered gaussian random variable with variance $1$, independent from all the former random variables, and set
\begin{equation*}
\phi_L^{\pm}=\pm\tilde{\xi}\frac{f}{\|\nabla f\|_2}+\tilde{\phi}_L.
\end{equation*}
Then, $\phi_L^{\pm}$ are random fields with the same law as $\phi_L$ but independent from $\xi$. Furthermore,
\begin{equation*}
\tilde{\phi}_L=\frac{\phi_L^-+\phi_L^+}{2}.
\end{equation*}
We now introduce a constant $A>0$ which we will fix later. The field $\phi_L$ will be positive on $D$ if the following three equations are satisfied.
\begin{align*}
\xi&>\|\nabla f\|_2 A\frac{1}{1-\eps}\\
\forall x\in D,\ \phi_L^-(x)&\leq A\\
\forall x\in D,\ \phi_L^+(x)&\leq A.\\
\end{align*}
Therefore, by independence of $\xi$ and $\phi_L^\pm$,
\begin{equation*}
\prob(\Omega_L^+)\geq \Big(1-2\prob\Big(\sup_D\phi_L>A\Big)\Big)\prob(\xi > \|\nabla f\|_2 A(1-\eps)^{-1}).
\end{equation*}
Choose $\delta>0$ and $A=\Big(\sqrt{\frac{2}{\pi}}+\delta\Big)\ln(\sqrt{L})$ for $L\geq L_0$. From Theorem \ref{supthm} we have
\begin{equation*}
\prob\Big(\sup_D\phi_L\geq\Big(\sqrt{\frac{2}{\pi}}+\delta\Big)\ln\big(\sqrt{L}\big)\Big)\rightarrow 0.
\end{equation*}
Moreover, from gaussian tail estimates (see equation (1.2.2) of \cite{adta_rfg}), for $L$ large enough,
\[\prob\Big(\xi > (1-\eps)^{-1}||\nabla f||_2\Big(\sqrt{\frac{2}{\pi}}+\delta\Big)\ln\big(\sqrt{L}\big)\Big)\]
is greater than
\[ \frac{(1-\eps)\exp\Big(-\frac{1}{2}(1-\eps)^{-2}\|\nabla f\|_2^2\Big(\sqrt{\frac{2}{\pi}}+\delta\Big)^2\ln^2\big(\sqrt{L}\big)\Big)}{2\|\nabla f\|_2\Big(\sqrt{\frac{2}{\pi}}+\delta\Big)\ln\big(\sqrt{L}\big)}.\]
Therefore,
\begin{align*}
\liminf_{L\rightarrow\infty}\frac{\ln\prob(\Omega_L^+)}{\ln^2\big(\sqrt{L}\big)}&\geq -\Big(\sqrt{\frac{2}{\pi}}+\delta\Big)^2\frac{1}{2}(1-\eps)^{-2}\|\nabla f\|_2^2\\
                                                                           &\geq -\Big(\sqrt{\frac{2}{\pi}}+\delta\Big)^2\frac{1}{2}(1-\eps)^{-2}(\|\nabla h\|_2+\eps)^2.
\end{align*}
Taking the infimum over $\delta>0$ and $\eps>0$ we obtain
\begin{equation*}
\liminf_{L\rightarrow\infty}\frac{\ln\prob(\Omega_L^+)}{\ln^2\big(\sqrt{L}\big)}\geq -\frac{2}{\pi}\frac{1}{2}\|\nabla h\|_2^2.
\end{equation*}
By taking the infimum of $\frac{1}{2}\|\nabla h\|_2^2$ over all $h$, we get
\begin{equation*}
\liminf_{L\rightarrow\infty}\frac{\ln\prob(\Omega_L^+)}{\ln^2\big(\sqrt{L}\big)}\geq -\frac{2}{\pi}cap(D).
\end{equation*}
\end{prf}

\subsection{The upper bound}

The aim of this section is to prove the following proposition. 
\begin{prop}\label{posub}
\begin{equation*}
\limsup_{L\rightarrow\infty}\frac{\ln\big(\Omega_L^+\big)}{\ln^2\big(\sqrt{L}\big)}\leq-\frac{2}{\pi}cap(D).
\end{equation*}
\end{prop}

Just as in section 3 of \cite{bdg_entrop}, we proceed by dichotomy with respect to the following event. Let $K\in\nats\setminus\{0\}$, $\eta>0$, $\alpha\in]0,1]$ and
\[A_{K,\eta,\alpha}=\Big\{Vol\Big[x\in D\ |\ \phi_{L^\alpha}(x) < \Big(\sqrt{\frac{2}{\pi}}-\eta\Big)\ln\big(\sqrt{L}\big)\Big]\leq\frac{K}{L^{\alpha/2}}\Big\}.\]
Proposition \ref{posub} is an immediate consequence of the two following results, to be compared with lemmas 9 and 10 of \cite{bdg_entrop}.

\begin{lm}\label{hubon}
Let $\eta>0$. For any integer $K>0$ and $\alpha\in]0,1]$,
\begin{equation*}
\limsup_{L\rightarrow\infty}\frac{\ln\big(\prob(A_{K,\eta,\alpha}\cap\Omega_L^+)\big)}{\ln^2\big(\sqrt{L}\big)}\leq-\Big(\sqrt{\frac{2}{\pi}}-\eta\Big)^2cap(D).
\end{equation*}
\end{lm}

\begin{lm}\label{huboff}
For any $\eta>0$ and $\kappa>0$ there exist $\alpha\in]0,1[$ and an integer $K>0$ such that
\begin{equation*}
\limsup_{L\rightarrow\infty}\frac{\ln\big(\prob(A_{K,\eta,\alpha}^c\cap\Omega_L^+)\big)}{\ln^2\big(\sqrt{L}\big)}\leq-\kappa.
\end{equation*}
\end{lm}

\begin{prf}[ of Lemma \ref{hubon}]
Fix $\eta,K>0$ and $\alpha\in]0,1]$. Choose $f\in C^\infty(\Sigma)$ positive ond $D$. For each $L>0$ let $\mathcal{E}_L$ be the random set and $X_L$ the real random variable defined as follows.
\begin{align*}
\mathcal{E}_L&=\Big\{x\in D\ |\ \phi_{L^\alpha}(x)\geq \Big(\sqrt{\frac{2}{\pi}}-\eta\Big)\ln\big(\sqrt{L}\big)\Big\}\\
X_L&=\int_D f(p)\phi_{L^\alpha}(p)|dV_g|(p).
\end{align*}
$X_L$ is a centered gaussian variable with variance
\[Var(X_L)=\int_\Sigma\int_\Sigma \mathbf{1}_D(p)f(p)G_{L^\alpha}(p,q)\mathbf{1}_D(q)f(q)|dV_g|(p)|dV_g|(q).\]
Therefore, according to equation \eqref{ref2},
\[ Var(X_L)\xrightarrow[L\to\infty]{} \sigma(\mathbf{1}_Df). \]
Moreover, on $A_{K,\eta,\alpha}\cap\Omega_L^+$, $Vol(D\setminus \mathcal{E}_L)\leq\frac{K}{L^{\alpha/2}}$ and $X_L\geq \Big(\sqrt{\frac{2}{\pi}}-\eta\Big)\ln\big(\sqrt{L}\big)\int_{\mathcal{E}_L}f$ so that
\begin{equation*}
X_L\geq \Big(\sqrt{\frac{2}{\pi}}-\eta\Big)\ln\big(\sqrt{L}\big)\Big(\int_Df-\frac{K||f||_\infty}{L^{\alpha/2}}\Big).
\end{equation*}

Therefore,
\[\prob(A_{K,\eta,\alpha}\cap\Omega_L^+)\leq \prob\Big(X_L\geq \Big(\int_Df-\frac{K||f||_\infty}{L^{\alpha/2}}\Big)\Big(\sqrt{\frac{2}{\pi}}-\eta\Big)\ln\big(\sqrt{L}\big)\Big).\]
By standard gaussian tail estimates (see again equation (1.2.2) of \cite{adta_rfg}), for $L$ large enough, the right hand side of the above inequality is smaller than
\[\frac{\sqrt{Var(X_L)}\exp\Big(-\frac{\Big(\int_Df-\frac{K||f||_\infty}{L^{\alpha/2}}\Big)^2}{2Var(X_L)}\Big(\sqrt{\frac{2}{\pi}}-\eta\Big)^2\ln^2\big(\sqrt{L}\big)\Big)}{\Big(\int_Df-\frac{K||f||_\infty}{L^{\alpha/2}}\Big)\Big(\sqrt{\frac{2}{\pi}}-\eta\Big)\ln\big(\sqrt{L}\big)}.\]
Thus,
\begin{equation*}
\limsup_{L\rightarrow\infty}\frac{\ln\big(\prob(A_{K,\eta,\alpha}\cap\Omega_L^+)\big)}{\ln^2\big(\sqrt{L}\big)}\leq-\Big(\sqrt{\frac{2}{\pi}}-\eta\Big)^2\frac{\Big(\int_Df\Big)^2}{2\sigma(\mathbf{1}_Df)}.
\end{equation*}
From Proposition \ref{capalt}, we complete the proof by taking the supremum over $f$ in the last inequality.
\end{prf}

To prove Lemma \ref{huboff} we will need the following technical result, which we prove at the end of the section. This result is analogous to Lemma 12 of \cite{bdg_entrop}.

\begin{lm}\label{fluct}
Let $\eta>0$ and $\alpha\in ]0,1]$. For each $\delta>0$, let $F_\delta$ be the event defined by
\begin{equation*}
F_\delta=\Big\{\sup_{d_g(p,q)\leq\delta L^{-\alpha/2}}|\phi_{L^\alpha}(p)-\phi_{L^\alpha}(q)|\geq (\eta/2)\ln\big(\sqrt{L}\big) \Big\}.
\end{equation*}
Then, for each $\kappa>0$ there is $\delta_0>0$ such that for each $\delta<\delta_0$,
\begin{equation*}
\prob\big(F_\delta\big) \leq \exp\big(-\kappa\ln^2\big(\sqrt{L}\big)\big).
\end{equation*}
\end{lm}

The proof itself goes roughly as follows. On the event $A_{K,\eta,\alpha}^c$ the field $\phi_{L^\alpha}$ takes low values on an abnormally large set $\mathcal{E}_L^c$. On this set we will consider $K$ small disks that will be so far apart from each other that, from Proposition \ref{indep}, the values of the field $\psi_{\alpha,L}$ on different disks will decorrelate. Now Lemma \ref{cast} tells us that for $\alpha$ small enough, on each of these disks, $\psi_{\alpha,L}$ will be likely to spike downwards and make the sum negative. Taking $K$ large enough will yield the desired inequality.
 
\begin{prf}[ of Lemma \ref{huboff}]
Let us begin by fixing $\eta>0$ and $\kappa>0$. We introduce constants $K>0$ and $\alpha\in]0,1[$ which we will fix later. As in the previous lemma, for every $L>0$ let
\begin{equation*}
\mathcal{E}_L=\Big\{x\in D\ |\ \phi_{L^\alpha}(x)\geq \Big(\sqrt{\frac{2}{\pi}}-\eta\Big)\ln\big(\sqrt{L}\big)\Big\}
\end{equation*}
By Lemma \ref{fluct} there is $\delta_0>0$ such that for any $\delta<\delta_0$,
\[\prob(\Omega^+_L\cap A_{K,\eta,\alpha}^c\cap F_\delta)\leq\prob(F_\delta)\leq\exp\Big(-\kappa\ln^2\big(\sqrt{L}\big)\Big).\]
We now provide an upper bound for $\prob(\Omega^+_L\cap A_{K,\eta,\alpha}^c\cap F_\delta^c)$. On the event $A_{K,\eta,\alpha}^c$, there is a constant $c>0$ such that for $L$ large enough, $ D \setminus \mathcal{E}_L$ contains at least $K$ points whose mutual distance and distance to $\mathcal{E}_L$ is at least $2c L^{-\alpha/4}$. When both $A_{K,\eta,\alpha}^c$ and $F_\delta^c$ are satisfied, for $L$ large enough, there is a collection $(D_j)_{j\in J}$ of $K$ disks of radius $\delta L^{-\alpha/2}$ on which $\phi_{L^\alpha}$ is smaller than $\Big(\sqrt{\frac{2}{\pi}}-(\eta/2)\Big)\ln\big(\sqrt{L}\big)$ such that for each $i,j\in J$ distinct,
\begin{equation*}
\inf_{p\in D_i,\ q\in D_j}d_g(p,q)\geq cL^{-\alpha/4}.
\end{equation*}
Recall that $\psi_{\alpha,L}$, defined in equation \eqref{ref4}, is independent from $\phi_{L^\alpha}$. Le $\tilde{\prob}$ be the conditional probability with respect to $\phi_{L^\alpha}$. On $A_{K,\eta,\alpha}^c\cap F_\delta^c$ consider the events
\begin{equation*}
T_{K,\eta,\alpha}^\pm=\Big\{\phi_{L^\alpha}\in A_{K,\eta,\alpha}^c\cap F_\delta^c\text{ and }\forall j\in J,\ \sup_{D_j}\pm\psi_{\alpha,L}\leq \Big(\sqrt{\frac{2}{\pi}}-(\eta/2)\Big)\ln\big(\sqrt{L}\big)\Big\}.
\end{equation*}
Since, $\psi_{\alpha,L}$ has the same law as $-\psi_{\alpha,L}$, the two events have the same probability. Moreover, since $\phi_L = \phi_{L^\alpha}+\psi_{\alpha,L}$, $\Omega^+_L\cap A_{K,\eta,\alpha}^c\cap F_\delta^c$ clearly implies $T_{K,\eta,\alpha}^{-}$. We will prove that for an adequate choice of $K$ and $\alpha$, and for $L$ large enough, on the event $A_{K,\eta,\alpha}^c\cap F_\delta^c$,
\begin{equation}\eqlab \label{ref10}
\tilde{\prob}\Big(T_{K,\eta,\alpha}^+\Big)\leq \exp\Big(-\kappa\ln^2\big(\sqrt{L}\big)\Big).
\end{equation}
Passing to expectations with respect to $\phi_{L^\alpha}$, this will imply that
\[\prob(\Omega^+_L\cap A_{K,\eta,\alpha}^c\cap F_\delta^c)\leq\prob\Big(T_{K,\eta,\alpha}^+\Big)\leq \exp\Big(-\kappa\ln^2\big(\sqrt{L}\big)\Big).\]
For each $j\in J$, consider $\iota_{j,\alpha,L}:V_{\floor{L^{(1-\alpha)/2}}}\rightarrow D_j$ such that for any distinct $x,y\in V_{\floor{L^{(1-\alpha)/2}}}$
\begin{equation*}
\frac{\delta}{8}\leq\frac{\sqrt{L}d_g(\iota_{\alpha,L}(x),\iota_{\alpha,L}(y))}{|x-y|}\leq \frac{\delta}{2}.
\end{equation*}
Here, as in Lemma \ref{cast}, $V_N$ is the box of sidelength $N$ in $\mathbb{Z}^2$. Now, for each $j\in J$, define $\psi_{j,\alpha,L}$ a random field on $V_{\floor{L^{(1-\alpha)/2}}}$ with the same law as $\psi_{\alpha,L}\circ\iota_{j,\alpha,L}$ such that the collection $(\psi_{j,\alpha,L})_{j\in J}$ is independent overall and of the previously defined fields. According to Lemma \ref{cast} and Proposition \ref{lcmax}, there is a constant $a>0$ depending only on $\eta$ and $\delta$ such that on the event $A_{K,\eta,\alpha}^c\cap F_\delta^c$,
\begin{equation}\eqlab \label{ref9}
\forall j\in J,\ \tilde{\prob}\Big(\sup\psi_{j,\alpha,L}\leq\Big(\sqrt{\frac{2}{\pi}}-(\eta/16)\Big)(1-\alpha)\ln\big(\sqrt{L}\big)\Big)\leq\exp\big(-a\ln^2\big(\sqrt{L}\big)\big).
\end{equation}
At this point, we fix $\alpha>0$ such that $\Big(\sqrt{\frac{2}{\pi}}-(\eta/16)\Big)(1-\alpha)\geq\Big(\sqrt{\frac{2}{\pi}}-(\eta/8)\Big)$. For each $j\in J$, define $V_j=\iota_{j,\alpha,L}(V_{\floor{L^{(1-\alpha)/2}}})$.  Take $\eps>0$. By Proposition \ref{indep} with $\beta=\alpha/2$, there is $L_0$ such that for any $L\geq L_0$, any distinct $i,j\in J$ and any $p\in V_i$, $q\in V_j$,
\begin{equation}\eqlab \label{ref6}
|\esp[\psi_{\alpha,L}(p)\psi_{\alpha,L}(q)]|\leq \eps^2.
\end{equation}
We now introduce, for each $j\in J$ and $p\in V_j$, a real random variable $\xi_p$, as well as an additional random variable $\xi$, all independent from previously introduced variables such that the family $(\xi,(\xi_p)_p)$, is independent and each variable is a centered gaussian with variance $1$. We will consider the following events

\begin{align*}
\Lambda_1&=\Big\{\forall j\in J, \sup_{p\in V_j}\psi_{\alpha,L}(p)\leq \Big(\sqrt{\frac{2}{\pi}}-(\eta/2)\Big)\ln\big(\sqrt{L}\big)\Big\}\\
\Lambda_2&=\Big\{\exists j\in J,\ p\in V_j\ |\ \eps\xi_p\geq (\eta/4)\ln\big(\sqrt{L}\big)\Big\}\\
\Lambda_3&=\Big\{\forall j\in J, \sup_{p\in V_j}\Big[\psi_{\alpha,L}(p)+\eps\xi_p\Big]\leq \Big(\sqrt{\frac{2}{\pi}}-(\eta/4)\Big)\ln\big(\sqrt{L}\big)\Big\}\\
\Lambda_4&=\Big\{\forall j\in J, \sup_{V_{\floor{L^{(1-\alpha)/2}}}}\Big[\psi_{j,\alpha,L} + \eps\xi\Big]\leq\Big(\sqrt{\frac{2}{\pi}}-(\eta/4)\Big)\ln\big(\sqrt{L}\big)\Big\}\\
\Lambda_5&=\Big\{\forall j\in J, \sup_{V_{\floor{L^{(1-\alpha)/2}}}}\psi_{j,\alpha,L} \leq\Big(\sqrt{\frac{2}{\pi}}-(\eta/8)\Big)\ln\big(\sqrt{L}\big)\Big\}\\
\Lambda_6&=\Big\{\eps\xi>(\eta/8)\ln\big(\sqrt{L}\big)\Big\}.
\end{align*}

We have the following inclusions

\begin{align*}
T_{K,\eta,\alpha}^+\subset\Lambda_1&\subset\Lambda_2\cup\Lambda_3\\
\Lambda_4&\subset\Lambda_5\cup\Lambda_6.
\end{align*}

Consider the two following gaussian fields defined over $\sqcup_{j\in J}V_j$. First $(\psi_{\alpha,L}(p)+\eps\xi_p)_p$, then $p\mapsto \psi_{j,\alpha,L}\circ\iota_{j,\alpha,L}^{-1}(p)+\eps\xi$ where $j\in J$ is such that $p\in V_j$. From equation \eqref{ref6}, these two fields satisfy the conditions for Lemma \ref{slep} so that
\begin{equation*}
\tilde{\prob}(\Lambda_3)\leq\tilde{\prob}(\Lambda_4).
\end{equation*}
Therefore,
\begin{equation*}\eqlab \label{ref7}
\tilde{\prob}\Big(T_{K,\eta,\alpha}^+\Big)\leq\tilde{\prob}(\Lambda_1)\leq\tilde{\prob}(\Lambda_2)+\tilde{\prob}(\Lambda_5)+\tilde{\prob}(\Lambda_6).
\end{equation*}
Firstly, by standard estimates on tails of gaussian variables, for $L$ large enough,
\begin{align*}\eqlab \label{ref8}
\tilde{\prob}(\Lambda_2)&\leq\frac{4K\eps\floor{L^{(1-\alpha)/2}}^2}{\eta\ln\big(\sqrt{L}\big)}\exp\Big(-\frac{\eta^2}{32\eps^2}\ln^2\big(\sqrt{L}\big)\Big)\\
\tilde{\prob}(\Lambda_6)&\leq\frac{8\eps}{\eta\ln\big(\sqrt{L}\big)}\exp\Big(-\frac{\eta^2}{128\eps^2}\ln^2\big(\sqrt{L}\big)\Big)
\end{align*}
Secondly, by independence of the $\psi_{j,\alpha,L}$ and by equation \eqref{ref9}
\begin{align*}
\tilde{\prob}(\Lambda_5)&= \prod_{j\in J}\tilde{\prob}\Big(\sup\psi_{j,\alpha,L}\leq\Big(\sqrt{\frac{2}{\pi}}-(\eta/8)\Big)\ln\big(\sqrt{L}\big)\Big).\\
                        &\leq \exp\big(-aK\ln^2\big(\sqrt{L}\big)\big).
\end{align*}
From equations \eqref{ref7} and \eqref{ref8}, taking $K$ large enough and $\eps$ small enough, we obtain inequality \ref{ref10} and the lemma is proved.
\end{prf}

We now prove Lemma \ref{fluct}. The proof is very similar to that of Proposition \ref{supup}.
\begin{prf}[ of Lemma \ref{fluct}]
Firstly, note that if $0<\delta<\delta'$ then $F_\delta\subset F_{\delta'}$. Therefore, for each $C>0$ it is enough to find $\delta>0$ such that $\prob(F_\delta)\leq e^{-C\ln\big(\sqrt{L}\big)^2}$. Choose $p,q\in\Sigma$ at distance smaller or equal to $\delta L^{-\alpha/2}$. There is a smooth path $\gamma$ in $\Sigma$, parametrized by arclength, such that $\gamma(0)=p$, $\gamma(2\delta L^{-\alpha/2})=q$ and for each $0\leq t\leq 2\delta L^{-\alpha/2}$, $|\gamma'(t)|=1$. Thus,
\begin{align*}
|\phi_{L^\alpha}(p)-\phi_{L^\alpha}(q)|&\leq \int_0^{2\delta L^{-\alpha/2}}|(\phi_{L^\alpha}\circ\gamma)'(t)|dt\\
                    &\leq 2\delta L^{-\alpha/2}\sup_{D_{L^\alpha}(p,2\delta)}|\nabla\phi_{L^\alpha}|.
\end{align*}

We now choose a local trivialisation the tangent bundle of $D_{L^\alpha}(p,2\delta)$ in which $\nabla\phi_{L^\alpha}$ has coordinates $(\nabla^1\phi_{L^\alpha},\nabla^2\phi_{L^\alpha})$. From Proposition \ref{kerprop}, there is a constant $C>0$ independent of $p$ such that for any such $q$ and for $j=1,2$,

\begin{equation*}
Var(\nabla^j_q\phi_{L^\alpha})\leq C L^{\alpha}.
\end{equation*}

With this information as well as the second inequality in Lemma \ref{sob}, we apply the BTIS inequality (Proposition \ref{bor}) to the random fields $(\nabla^j_q\phi_{L^\alpha}(v_q))_{q\in D_{L^\alpha}(p,2\delta)}$ for $j=1,2$ and deduce there is a constant $C>0$ such that for $L$ large enough,

\begin{equation*}
\prob\Big(\sup_{D_{L^\alpha}(p,2\delta)}|\nabla\phi_{L^\alpha}|\geq \frac{\eta L^{\alpha/2}}{4\delta}\ln\big(\sqrt{L}\big)\Big)\leq \exp\Big(-\frac{C\eta^2}{\delta^2}\ln^2\big(\sqrt{L}\big)\Big).
\end{equation*}

By the triangle inequality,

\begin{equation}\eqlab \label{ref5}
\prob\Big(\sup_{q_1,q_2\in D_{L^\alpha}(p,\delta)}|\phi_{L^\alpha}(q_1)-\phi_{L^\alpha}(q_2)|\geq (\eta/2)\ln\big(\sqrt{L}\big)\Big)\leq \exp\Big(-\frac{C\eta^2}{\delta^2}\ln^2\big(\sqrt{L}\big)\Big).
\end{equation}

There exists $C>0$ such that for each $L>0$, there is a covering of $\Sigma$ by a collection $(D_j)_{j\in J}$ of at most $C L^\alpha \delta^{-2}$ disks of radius $\delta L^{-\alpha/2}$. For each $j$, inequality \eqref{ref5} applies on $D_j$. Therefore,

\begin{align*}
\prob\Big(\sup_{d_g(p,q)\leq \delta L^{-\alpha/2}}|\phi_{L^\alpha}(p)-\phi_{L^\alpha}(q)|\geq (\eta/2)\ln\big(\sqrt{L}\big)\Big)&\leq C L^\alpha \exp\Big(-\frac{C\eta^2}{\delta^2}\ln^2\big(\sqrt{L}\big)\Big)\\
                                                                                                                &\leq \exp\Big(-\frac{2C\eta^2}{\delta^2}\ln^2\big(\sqrt{L}\big)\Big)
\end{align*}
for $L$ large enough. This inequality ends the proof of the lemma.
\end{prf}

\section{The covariance function}\label{sect_ana}

The aim of this section is to prove Theorem \ref{kerthm}, Proposition \ref{kerprop} and Proposition \ref{indep}. For each $L>0$, let $E_L$ be the Schwartz kernel of the orthogonal projector in $L^2(\Sigma)$ onto $U_L$. Then,
\begin{equation}\label{kereq2}
E_L(p,q)=\sum_{0<\lambda_n\leq L}\psi_n(p)\psi_n(q).
\end{equation}
The asymptotic behavior of this kernel as $L\rightarrow\infty$ has been studied extensively by H\"ormander in \cite{ho_sfeo} (see also \cite{bin_der} and \cite{gawe_betti}). In order to prove Theorem \ref{kerthm} we will express $G_L$ in terms of $E_L$ and use the aforementioned results to extract an explicit formula for $G_L$.

\subsection{Preliminary results}

We begin with the following proposition, which establishes the link between $G_L$ and $E_L$.

\begin{prop}\label{link}
There is a function $R\in C^\infty(M\times M, \reals)$ such that for each $L>0$,
\begin{equation*}
G_L = \frac{E_L}{L}+\int_1^L\frac{E_\lambda}{\lambda^2}d\lambda + R.
\end{equation*}
\end{prop}

\begin{prf}
For any $p,q\in \Sigma$, the functions $L\mapsto E_L(p,q)$ and $L\mapsto G_L(p,q)$ define distributions over $]0,+\infty[$. In what follows, $\partial_L$ will mean differentiation in the sense of distributions. First of all, for any $p,q\in\Sigma$

\begin{align*}
\partial_LE_L(p,q)&=\sum_{0<\lambda_n}\psi_n(p)\psi_n(q)\delta_{\lambda_n}(L)\\
\partial_LG_L(p,q)&=\sum_{0<\lambda_n}\frac{1}{\lambda_n}\psi_n(p)\psi_n(q)\delta_{\lambda_n}(L)=\sum_{0<\lambda_n}\frac{1}{L}\psi_n(p)\psi_n(q)\delta_{\lambda_n}(L)=\frac{1}{L}\partial_LE_L(p,q).
\end{align*}
Consequently $\partial_L\Big(G_L-\frac{E_L}{L}+\int_1^L\frac{E_\lambda}{\lambda^2}d\lambda\Big)=0$. Therefore
\begin{equation*}
R=G_L-\frac{E_L}{L}+\int_1^L\frac{E_\lambda}{\lambda^2}d\lambda
\end{equation*}
is independent of $L$. The right hand side is a linear combination of functions $(p,q)\mapsto\psi_n(p)\psi_n(q)$ so it belongs to $C^\infty(M\times M)$.
\end{prf}

Now, we use Theorem 5.1 of \cite{ho_sfeo} to obtain an explicit description of the integral term in the equation of the previous proposition. Let us fix $p_0\in \Sigma$. According to Theorem 5.1 of \cite{ho_sfeo} there is an open neighborhood $U$ of $p$ in $\Sigma$ with a chart $\phi:U\rightarrow\phi(U)$ such that $(d_{p_0}\phi^*)^{-1}$ is an isometry from $T^*_{p_0}\Sigma$ with the metric induced by $g$ to $\reals^n$ with the euclidian metric, as well as a real valued function $\theta\in C^\infty(U\times T^*U)$ satisfying the phase condition (see Definition 2.3 of \cite{ho_sfeo}) and a constant $C>0$ such that for each $p,q\in U$ and each $L>0$,
\begin{equation}\label{ho_eq}
\Big|E_L(p,q) - \frac{1}{(2\pi)^2}\int_{|\xi|^2\leq L} e^{i\theta(p,q,\xi)}d_q\eta(\xi)\Big|\leq C\sqrt{L}
\end{equation}
where $d_q\eta$ is the measure associated to the metric induced by $g$ on $T^*U$. Let $\phi_*\theta \in C^\infty(\phi(U)\times\phi(U)\times \reals^2)$ be defined by
\[\forall x,y\in\phi(U),\ w\in\reals^2,\ \phi_*\theta(x,y,w)=\theta(\phi^{-1}(x),\phi^{-1}(y),d_{\phi^{-1}(y)}\phi^*w).  \]
Here $d_{\phi^{-1}(y)}\phi^*$ is the adjoint of the differential of $\phi$ at $\phi^{-1}(y)$.\\

Theorem 5.1 of \cite{ho_sfeo} provides the following information concerning $\theta$.
\begin{enumerate}
\item For each $x,y\in \phi(U)$ and $w\in\reals^2$ such that $\langle x-y,w\rangle = 0$, $\phi_*\theta(x,y,w)=0$.
\item For each $y\in \phi(U)$ and $w\in\reals^2$, $\partial_x(\phi_*\theta)(x,y,w)|_{x=y}=w$.
\end{enumerate}
In particular, these equations have the following consequence. Choose $y\in\phi(U)$, $t\geq 0$, $v\in S^1$, $\alpha\in \nats^2$ a multiindex and $w\in\reals^2$. Then, by the Taylor-Young estimate applied to $\partial_w^\alpha\phi_*\theta(y+\lambda v,y,w)$ with respect $\lambda$, for $\lambda$ small enough,
\begin{equation}\label{ho_ph}
\partial_w^\alpha\phi_*\theta(y+\lambda v,y,w)=\lambda\partial_w^\alpha\langle v,w\rangle + O(\lambda^2)
\end{equation}
where the $O(\lambda^2)$ is uniform when $y$ and $w$ are restricted to any compact set.\\

Before we proceed any further, let us introduce some notation. For each $q\in\Sigma$, let $S_q$ be the unit circle in $(T_q^*\Sigma,g_q)$ and $d_q\nu$ the measure induced by the restriction of $g_q$ to $S_q$. Also, for each $y\in\phi(U)$, let $\tilde{S}_y$ be
\begin{equation}\eqlab \label{ref23}
\tilde{S}_y=\{w\in\reals^2 |\ d_{\phi^{-1}(y)}\phi^*w\in S_{\phi^{-1}(y)}\}.
\end{equation}
\begin{lm}\label{fubi}
For each $t>0$ and $p,q\in U$, let
\begin{equation}\label{jdef}
J(p,q,t)=\int_{S_q}e^{i\theta(p,q,t\omega)}d_q\nu(\omega).
\end{equation}
Then, there is a constant $C>0$ such that for each $L>0$ and $p,q\in U$,
\begin{equation*}
\Big|\int_1^LE_\lambda(p,q)\lambda^{-2}d\lambda-\frac{1}{4\pi^2}\Big(\int_1^{\sqrt{L}}J(p,q,t)t^{-1}dt-L^{-1}\int_1^{\sqrt{L}}J(p,q,t)tdt\Big)\Big|\leq C.
\end{equation*}
\end{lm}

\begin{prf}
First of all, from equation (\ref{ho_eq}) there is $C>0$ such that for all $p,q$ and $L$,

\begin{align*}
\Big| \int_1^LE_\lambda(p,q)\lambda^{-2}d\lambda-\frac{1}{4\pi^2}\int_1^L\lambda^{-2}\int_{|\xi|^2\leq \lambda}e^{i\theta(p,q,\xi)}d_p\eta(\xi)d\lambda\Big| &\leq C\int_1^L\lambda^{-3/2}d\lambda\\
                                                                                                                                                              &\leq 2C.\eqlab \label{fubi_1}
\end{align*}

Now, applying the polar change of coordinates $(t,\omega)\mapsto t\omega=\xi$,

\[ \int_{|\xi|^2\leq \lambda}e^{i\theta(p,q,\xi)}d_p\eta(\xi) = \int_0^{\sqrt{\lambda}} J(p,q,t)tdt. \]

Next, we apply the change of variables $u=\sqrt{\lambda}$.

\[ \int_1^L\lambda^{-2} \int_0^{\sqrt{\lambda}} J(p,q,t)tdtd\lambda = \int_1^{\sqrt{L}}\int_0^uJ(p,q,t)tdt 2u^{-3}du. \]

We split the inner integral in two $\int_0^u=\int_0^1+\int_1^u$. Note that the integral from $0$ to $1$ has lost any dependence on $u$ or $L$. Since $\int_0^\infty u^{-3/2}du<\infty$ that term is bounded. Therefore applying Fubini's theorem,

\begin{align*}
\frac{1}{4\pi^2}\int_1^L\lambda^{-2}\int_0^{\sqrt{\lambda}}J(p,q,t)tdtd\lambda&=\frac{1}{4\pi^2}\int_1^{\sqrt{L}}\int_1^uJ(p,q,t)tdt2u^{-3}du\\
                                                                         &=\frac{1}{4\pi^2}\int_1^{\sqrt{L}}J(p,q,t)t\int_t^{\sqrt{L}}2u^{-3}dudt\\
                                                                         &=\frac{1}{4\pi^2}\Big(\int_1^{\sqrt{L}}J(p,q,t)t^{-1}dt-L^{-1}\int_1^{\sqrt{L}}J(p,q,t)tdt\Big).
\end{align*}
Together with equation (\ref{fubi_1}) this implies that
\begin{equation*}
\int_1^LE_\lambda(p,q)\lambda^{-2}d\lambda-\frac{1}{4\pi^2}\Big(\int_1^{\sqrt{L}}J(p,q,t)t^{-1}dt-L^{-1}\int_1^{\sqrt{L}}J(p,q,t)tdt\Big)
\end{equation*}
is bounded as required.
\end{prf}

To proceed any further, we need to control the behavior of $J(p,q,t)$ when $t\rightarrow\infty$. We will use the stationary phase method to prove the following proposition.
 
\begin{prop}\label{jprop}
There exist $V\subset U$ an open neighborhood of $p_0$ and a constant $C>0$, such that for all $p,q\in V$ and $t\in[0,+\infty[$,
\begin{equation}\label{jbound}
|J(p,q,t)|\leq \frac{C}{\sqrt{1+d_g(p,q)t}}.
\end{equation} 
\end{prop}

The function $J$ is an oscillatory integral over the circle. To obtain the bound for large $t$ uniformly with respect to $p,q$ distinct, we should apply the stationary phase method to the phase $w\mapsto \frac{\phi_*\theta(x,y,w)}{|x-y|}$ with parameter $|x-y|t$ where $\phi$ is the chart appearing just above equation (\ref{ho_eq}) and $\phi_*\theta$ is defined just below it. In order to obtain these bounds we will first apply the adequate change of variables in order to compactify the space of definition near the diagonal.

\begin{prf}[ of Proposition \ref{jprop}]
Let $K\subset \phi(U)$ be a compact neighborhood of $\phi(p)$. Since $K$ is compact, there exists a constant $\alpha>0$ such that for each $x\in K$, $\overline{D(x,\alpha)}\subset\phi(U)$. Here $D(x,\alpha)$ is the open disk centered at $x$ and of radius $\alpha$. Let us define the following sets.
\begin{align*}
A&=\{(x,y,w)\ |\ x,y\in K,\ x\neq y,\ |x-y|<\alpha,\ w\in \tilde{S}_y\}\\
B&= S^1\times]0,\alpha[\times\{(y,w)\in K\times\reals^2\ |\ y\in K,\ w\in\tilde{S}_y \}.
\end{align*}
where $\tilde{S}_y$ is defined in equation \eqref{ref23}. The following map is a diffeomorphism
\begin{align*}
f: A&\rightarrow B\\
  (x,y,w)&\mapsto\Big(\frac{x-y}{|x-y|},|x-y|,y,w\Big).
\end{align*}
Set
\begin{align*}
\psi: A&\rightarrow \reals\\
     (x,y,w)&\rightarrow\frac{\phi_*\theta(x,y,w)}{|x-y|}.
\end{align*}
Now, applying equation (\ref{ho_ph}) with $\alpha=0$ we deduce that for any $(v,\lambda,y,w)\in B$,

\[ \psi\circ f^{-1}(v,\lambda,y,w) = \langle v,w\rangle + O(\lambda) \]

where $O(\lambda)$ is uniform with respect to $v,y,w$. Thus, $\Psi:=\psi\circ f^{-1}$ extends by continuity to $\overline{B}$ so that for all $y\in K$, $w\in \tilde{S}_y$ and $v\in S^1$, $\Psi(v,0,y,w)=\langle v,w\rangle$. Equation (\ref{ho_ph}) also shows that the map

\[\lambda \mapsto \Big((v,y,w)\mapsto \Psi(v,\lambda,y,w)\Big) \]

is continuous from $[0,\alpha]$ to the space of continuous functions on $S^1\times\{(y,w)\in K\times\reals^2\ |\ w\in \tilde{S}_y\}$ which are $C^\infty$ with respect to $w$. Let us momentarily fix $v$ and set $y=\phi(p_0)$. The curve $\tilde{S}_y$ is a circle since $d_{p_0}\phi^{-1}$ is an isometry. Hence, the map $w\mapsto\Psi(v,0,y,w)=\langle v,w\rangle$ defined over $\tilde{S}_y$ is a Morse function with two critical points. Thus, for $\lambda$ small enough and $y$ close enough to $\phi(p_0)$, $w\mapsto\Psi(v,\lambda,y,w)$ is also a Morse function with two critical points. Moreover, it depends $C^4$-continuously on $v$, $t$, and $y$. Therefore there exists $\beta>0$ and a constant $C>0$ such that for any $x,y\in D(\phi(p_0),\beta)$ distinct we have the following.
\begin{itemize}
\item $\psi_{x,y}:w\mapsto\psi(x,y,w)$ is a Morse function with two critical points $z_1(x,y)$ and $z_2(x,y)$.
\item For $j=1,2$, $\big|det(d^2_{z_j(x,y)}\psi_{x,y})\big|>C^{-1}$ (Here $d^2_{z_j(x,y)}\psi_{x,y}$ is the hessian of $\psi_{x,y}$ which is well defined since $z_j(x,y)$ are critical points, see for instance the definition following Lemma 1.6 of \cite{nic_morse}).
\item $\|\psi_{x,y}\|_{C^4}\leq C$.
\end{itemize}
Consequently, $u=1$ and $f=\psi_{x,y}$ satisfy the hypotheses of Theorem 7.7.5 of \cite{ho_apdo} with $k=1$, from which we deduce that there is $C>0$ such that for these same $x$, $y$ and $\lambda$,
\[J(\phi^{-1}(x),\phi^{-1}(y),\lambda)\leq \frac{C}{\sqrt{|x-y|\lambda}}.\]
Moreover, for all $p$, $q$ and $t$, $|J(p,q,t)|\leq \int_{S_q}  1 d_q\nu = 2\pi$. Let $V=\phi^{-1}(D(\phi(p_0,\beta))$. From the two last inequalities, there is a constant $C>0$ such that for any $p,q\in V$ and $t\geq 0$,
\[J(p,q,t)\leq \frac{C}{\sqrt{1+d_g(p,q)t}}.\]
\end{prf}

\subsection{Proof of Theorem \ref{kerthm}}

We now use the results and notations of the previous section to prove Theorem \ref{kerthm}. We will estimate $G_L$ when $p$ and $q$ are in a neighborhood of $p_0$ and use the compactness of $\Sigma$ to make the result global.. Firstly, from Proposition \ref{link}, $G_L = \frac{E_L}{L}+\int_1^L\frac{E_\lambda}{\lambda^2}d\lambda + O(1)$. From equation (\ref{ho_eq}), $L^{-1}E_L$ is bounded. Now, from Lemma \ref{fubi}, there is $C>0$ such that for each $p,q\in U$
\[\Big|\int_1^LE_\lambda(p,q)\lambda^{-2}d\lambda-\frac{1}{4\pi^2}\Big(\int_1^{\sqrt{L}}J(p,q,t)t^{-1}dt-L^{-1}\int_1^{\sqrt{L}}J(p,q,t)tdt\Big)\Big|\leq C.\]
Let $V\subset U$ be as in Proposition \ref{jprop}. From equation (\ref{jbound}), there is a constant $C>0$ such that for all $p,q\in V$,
\[ \Big|L^{-1}\int_1^{\sqrt{L}}J(p,q,t)tdt\Big|\leq C. \]
Note that, from equation (\ref{ho_ph}), $\theta(p,q,0)=0$ so that $J(p,q,0)=2\pi$. Choose $p,q\in U$ distinct and set $r=d_g(p,q)$. Then, applying the change of variables $a=rt$,
\begin{align*}
\int_1^{\sqrt{L}}J(p,q,t)t^{-1}dt&=\int_r^{r\sqrt{L}}J(p,q,r^{-1}a)a^{-1}da\\
                              &=2\pi\int_r^{r\sqrt{L}}\frac{\mathbf{1}_{a\leq 1}}{a}da + \int_r^{r\sqrt{L}}\frac{J(p,q,r^{-1}a)-2\pi\mathbf{1}_{a\leq 1}}{a}da
\end{align*}
where $\mathbf{1}_{a\leq 1}$ equals $1$ if $a\leq 1$ and $0$ otherwise. Recall the notation introduced in equation \eqref{ref23}. Now, from equation (\ref{ho_ph}), for $0 < a\leq r$ and $x,y\in \phi(V)$, uniformly for $w\in \tilde{S}_y$,
\[\phi_*\theta(x,y,(a/r)w)=\frac{\langle x-y,w\rangle a}{r} + O(|x-y|^2a/r)\]
so $f:a\mapsto\frac{J(p,q,r^{-1}a)-2\pi\mathbf{1}_{a\leq 1}}{a}$ admits a continuous extension at $a=0$ equal to
\[ \int_{S_q} \frac{i\langle \phi(p)-\phi(q),(d_q\phi^*)^{-1}\omega\rangle}{r}d_q\nu(\omega).\]
Let $W\subset V$ be an open neighboorhood of $p_0$ on which $\phi$ is bi-lipschitz and recall that $r=d_g(p,q)$. Then, there exists $C>0$ independent of $L$ such that for all $p,q\in W$, $|\phi(p)-\phi(q)|\leq Cr$ and such that for all $q\in W$ and $\omega\in S_q$, $\|(d_q\phi^*)^{-1}\omega\|_{eucl}\leq C$. (Here, $\|\ \|_{eucl}$ denotes the euclidean norm on $\reals^2$.) Thus, the aforementioned continuous extension of $f$ is uniformly bounded for $p,q\in W$. Moreover, by equation (\ref{jbound}), $f$ is $O(a^{-3/2})$ uniformly with respect to $r$ when $a\rightarrow\infty$. Therefore it is integrable with uniform bounds over $p$ and $q$. Consequently, for any distinct $p,q\in W$ such that $d_g(p,q)<1$,
\begin{align*}
\int_1^{\sqrt{L}}\lambda^{-2}E_\lambda d\lambda &= \frac{1}{2\pi} \int_{d_g(p,q)}^{d_g(p,q)\sqrt{L}}\frac{\mathbf{1}_{a\leq 1}}{a}da + O(1)\\
                                           &= \frac{1}{2\pi} \Big(\ln\big(\min(1,\sqrt{L}d_g(p,q))\big)-\ln\big(\min(1,d_g(p,q))\big)\Big) + O(1)\\
                                           &= \frac{1}{2\pi}\Big(\ln\big(\sqrt{L}\big)-\ln_{+}\big(\sqrt{L}d_g(p,q)\big)\Big)+O(1).
\end{align*}
Here, the bounds implied by the $O$'s are uniform with respect to $p$ and $q$. For the last equality we use the fact that $\ln_{+}\big(d_g(p,q)\big)$ is bounded for $p,q\in W$. The case $p=q$ follows by continuity. Moreover, by compactness, one can cover $\Sigma$ with a finite number of such $W$'s so that there is a constant $\eps>0$ independent of $L$ for which the constants are uniform with respect to any $p,q$ such that $d_g(p,q)<\eps$.

\subsection{Proof of Propositions \ref{kerprop} and \ref{indep}}

\begin{prf}[ of Proposition \ref{kerprop}]
From Theorem 2.3 of \cite{gawe_betti} with $M=\Sigma$ and $P=\Delta$, there is a constant $C>0$ such that for all $p\in \Sigma$ and $L\geq 0$,

\[ |(Q_1\otimes Q_2)E_L(p,p)|\leq C\big(1+L^{1+d/2}\big) \]

Therefore, from Proposition \ref{link} for all $p\in\Sigma$ and $L>0$,

\begin{align*}
|(Q_1\otimes Q_2)G_L(p,p)|&\leq L^{-1}|(Q_1\otimes Q_2)E_L(p,p)| +\\
                          &+ \int_1^L\lambda^{-2}|(Q_1\otimes Q_2)E_\lambda(p,p)|d\lambda +|(Q_1\otimes Q_2)R(p,p)|\\
                         &\leq C\Big(1+L^{d/2}+\int_1^L \lambda^{-1+d/2}d\lambda\Big)\\
                         &\leq C'\big(1+L^{d/2}\big).
\end{align*}
\end{prf}

In order to prove Proposition \ref{indep}, we will need the following technical result, which we deduce from Theorem 5.1 of \cite{ho_sfeo}.

\begin{lm}\label{boundaway}
For each $\beta>0$ and $\delta>0$, there is $C>0$ such that for each $L,\lambda>0$ and $p,q\in\Sigma$ such that $d_g(p,q)\geq\delta L^{-\beta/2}$,

\[\big|E_\lambda(p,q)\big|\leq C L^{\beta/4}\lambda^{3/4}.\]
\end{lm}

\begin{prf}
First of all, by Theorem 5.1 of \cite{ho_sfeo}, there is $\eps>0$ and $C$ such that for each $p_0\in\Sigma$ and $q\in\Sigma$,
\begin{enumerate}
\item If $d_g(p_0,q)\geq \eps$ then, for all $\lambda>0$, $\big|E_\lambda(p_0,q)|\leq C\sqrt{\lambda}$ (see equation (5.4) of \cite{ho_sfeo}).
\item There is a chart $(\phi,U)$ such that $D(p_0,\eps)\subset U$, around $p_0$ as well as a function $\theta$ such that equations (\ref{ho_ph}) and (\ref{ho_eq}) hold with the same constant $C>0$.
\end{enumerate}
From 1. we need only deal with the case where $p_0$ and $q$ are $\eps$-close. Now, by a polar coordinate change, and with the same notations as in Lemma \ref{fubi}, for all $\lambda>0$,

\[\int_{|\xi|^2\leq \lambda}e^{i\theta(p_0,q,\xi)}d_q\eta(\xi)=\int_0^{\sqrt{\lambda}}J(p_0,q,t)tdt.\]

By equation (\ref{jbound}) there is $C>0$ independent of $p_0$ and $q$ such that for all $\lambda>0$ and all $q\neq p_0$,
\[\Big|\int_0^{\sqrt{\lambda}}J(p_0,q,t)tdt\Big|\leq C\frac{\lambda^{3/4}}{\sqrt{d_g(p_0,q)}}.\]
This concludes the proof.
\end{prf}

\begin{prf}[ of Proposition \ref{indep}]
From Proposition \ref{link},
\begin{equation*}
G_{L^{\alpha},L}=G_L-G_{L^\alpha}=\frac{E_L}{L}-\frac{E_{L^\alpha}}{L^{\alpha}}+\int_{L^\alpha}^L\frac{E_\lambda}{\lambda^2}d\lambda.
\end{equation*}
Now from Lemma \ref{boundaway} there is a constant $C>0$ such that for any $L>0$ and any $p,q\in\Sigma$ such that $d_g(p,q)\geq \delta L^{-\beta/2}$,
\begin{align*}
\big|L^{-1}E_L(p,q)\big|&\leq CL^{(\beta-1)/4}\\
\big|L^{-\alpha}E_{L^\alpha}(p,q)\big|&\leq CL^{(\beta-\alpha)/4}\\
\Big|\int_{L^\alpha}^L\frac{E_\lambda}{\lambda^2}d\lambda\Big|&\leq CL^{\beta/4}\int_{L^\alpha}^L\lambda^{-5/4}d\lambda\leq 2CL^{(\beta-\alpha)/4}
\end{align*}
and each term tends to $0$ since $0<\beta<\alpha<1$.
\end{prf}

\section{The case of a surface with boundary}\label{sect_bd}

In this paper, we studied a random zero mean field over a compact surface $\Sigma$ with $\partial\Sigma=\emptyset$. The bulk of our results remain valid in some sense if $\Sigma$ has a smooth boundary $\partial\Sigma\neq\emptyset$. The zero mean condition is then replaced by a Dirichlet zero boundary condition. In this section we state definitions and results for this setting and explain the minor modifications needed in each proof. Be advised that \textbf{we use the same notations as before to denote slightly different objects}.

\subsection{Definitions}

Let $(\Sigma,g)$ be a smooth surface with smooth boundary $\partial\Sigma$ equipped with a riemannian metric. We will denote the interior of $\Sigma$ by $\mathring{\Sigma}=\Sigma\setminus\partial\Sigma$. If $E(\Sigma)$ is a topological vector-space of functions on $\Sigma$ of which the $C^\infty(\Sigma)$ is a dense subspace, $E_0(\Sigma)$ will denote the closure in $E(\Sigma)$ of smooth functions with compact support in $\mathring{\Sigma}$. We begin by defining the CGFF on $\Sigma$. As in the introduction, the bilinear form $\langle u,v \rangle_\nabla:=\int_\Sigma g(\nabla u,\nabla v) |dV_g|$ defines a scalar product on $H^1_0(\Sigma)$. From Theorem 4.43 of \cite{gahulaf} there exist $(\psi_n)_{n\geq 1}\in C^\infty_0(\Sigma)^\nats$ and $(\lambda_n)_{n\geq 1}\in\reals^\nats$ such that $0<\lambda_1\leq\lambda_2\dots,\lambda_n\xrightarrow[n\to\infty]{}+\infty$, such that $(\psi_n)_{n\geq 1}$ is a Hilbert basis for $L^2_0(\Sigma)$ and such that for each $n\in\nats$, $\Delta \psi_n=\lambda_n\psi_n$. Then $\Big(\frac{1}{\sqrt{\lambda_n}}\psi_n\Big)_{n\geq 1}$ is a Hilbert basis of $(H^1_0,\langle\ ,\ \rangle_\nabla)$. For each $L>0$, let $(U_L,\langle\ ,\ \rangle_\nabla)$ be the subspace of $(H^1_0,\langle\ ,\ \rangle_\nabla)$ spanned by the functions $\frac{1}{\sqrt{\lambda_n}}\psi_n$ such that $\lambda_n\leq L$. Let $(\xi_n)_{n\in\nats}$ be an i.i.d sequence of centered gaussians with variance one. Then, for each $L>0$ we define the CGFF on $\Sigma$ as
\begin{equation}
\phi_L=\sum_{\lambda_n\leq L}\frac{\xi_n}{\sqrt{\lambda_n}}\psi_n.
\end{equation}
For each $L>0$ and $p,q\in\Sigma$, let $G_L(p,q)=\esp[\phi_L(p)\phi_L(q)]$.\\

Let $D$ be an open subset of $\mathring{\Sigma}$ at positive distance from $\partial\Sigma$. Then, we define the capacity of $D$ relative to $\Sigma$ as the infimum of $\frac{1}{2}\|\nabla h\|^2$ over all $h\in C^\infty_0(\mathring{\Sigma})$ and denote it by $cap_\Sigma(D)$ or $cap(D)$. If $D$ is non-empty and has smooth boundary, we will see that $cap(D)>0$. Note that in this setting, the capacity is conformally invariant. More precisely, consider $(\Sigma_1,g_1)$ and $(\Sigma_2,g_2)$ are two riemannian surfaces with smooth boundary such that there exists a conformal isomorphism, $f:\Sigma_1\rightarrow\Sigma_2$. Let $D_1\subset \Sigma_1$ and $D_2=f(D_1)$. We claim that $cap_{\Sigma_1}(D_1)=cap_{\Sigma_2}(D_2)$. Indeed, for each $h\in C^\infty_0(\Sigma_2)$ such that $h\geq 1$ on $D_2$, we have $h\circ f\in C^\infty_0(\Sigma_1)$ and $h\circ f\geq 1$ on $D_1$. The same is true for $f^{-1}$ if we exchange $1$ and $2$ so $f$ defines a bijection between the sets whose infimum define the capacities. Moreover, since $f$ is a conformal map, for any $p\in\Sigma_1$ and $v\in T_p\Sigma_1$, $g_2(d_pfv,d_pfv)=|det(d_pf)|g_1(v,v)$ where $det(d_pf)$ is the determinant of the matrix of $d_pf$ written in orthonormal bases. It follows that from the change of variables formula that for any $h\in C^{\infty}_0(\Sigma_2)$, $\|h\|_\nabla=\|h\circ f\|_\nabla$. This proves the claim.

\subsection{Main results}

\begin{thm}\label{mainthmz}
Let $(\Sigma,g)$ be a smooth compact surface with smooth non-empty boundary and let $(\phi_L)_{L>0}$ be the CGFF on $(\Sigma,g)$. Let $D$ be a non-empty open subset of $\mathring{\Sigma}$ with smooth boundary and positive distance to $\partial\Sigma$. Then,
\[ \lim_{L\rightarrow\infty}\frac{\ln\big(\prob(\forall x\in D, \phi_L(x)>0)\big)}{\ln^2\big(\sqrt{L}\big)}=-\frac{2}{\pi}cap(D). \]
\end{thm}

The only change in section \ref{sect_hole}, which contains the heart of the argument, is in the proof of Proposition \ref{capalt}. Though the statement is the same, $\Sigma$ has a boundary and the definition of $cap(D)$ has changed. For a proof in this case we refer the reader to Lemma 2.1 of \cite{bd_crit} where one should simply replace $V$ by $\Sigma$ and $\langle f,\phi\rangle_V - \frac{1}{2}\langle f, K_V f\rangle_V$ by $\frac{1}{2\sigma(\mathbf{1}_Df)}\Big(\int_D f\Big)^2$.

\begin{thm}\label{kerthmz}
Let $(\Sigma,g)$ be a compact riemannian surface with smooth non-empty boundary and let $(\phi_L)_{L>0}$ be the CGFF on $(\Sigma,g)$. 
For each $L> 0$ and $p,q\in\Sigma$, let
\[ G_L(p,q)=\esp[\phi_L(p)\phi_L(q)]. \]
Then, there exists $\eps>0$ such that for each $p,q\in \mathring{\Sigma}$ satisfying $d_g(p,q)\leq \eps$ and for each $L>0$,
\[ G_L(p,q)=\frac{1}{2\pi}\Big(\ln\big(\sqrt{L}\big)-\ln_+\big(\sqrt{L}d_g(p,q)\big)\Big) + \rho_L(p,q) \]
where, for any compact subset $K\subset\mathring{\Sigma}$, $\rho_L(p,q)$ is uniformly bounded for $L>0$ and $p,q\in K$.
\end{thm}

For the proof of Theorem \ref{kerthmz} we fix a compact subset $K\subset\mathring{\Sigma}$ and proceed as in the original proof but replace the statements ``uniform with respect to $p,q\in\Sigma$'' by ``uniform with respect to $p,q\in K$''. This works because Theorem 5.1 of \cite{ho_sfeo} is valid on non-compact manifolds except the bound on the remainder term is only uniform on compact sets.

\begin{thm}\label{supthmz}
Let $(\Sigma,g)$ be a smooth compact riemannian surface with non-empty boundary and let $(\phi_L)_{L>0}$ be the CGFF on $(\Sigma,g)$. Let $D$ be a non-empty open subset of $\mathring{\Sigma}$ at positive distance of $\partial\Sigma$ and $K$ be a compact subset of $\mathring{\Sigma}$. Then for each $\eta >0$,
\begin{equation*}
\limsup_{L\rightarrow\infty}\frac{\ln\Big(\prob\Big(\sup_K\phi_L>\Big(\sqrt{\frac{2}{\pi}}+\eta\Big)\ln\big(\sqrt{L}\big)\Big)\Big)}{\ln\big(\sqrt{L}\big)}\leq -2\sqrt{2\pi}\eta+O(\eta^2)
\end{equation*}
and there exists $a>0$ such that for $L$ large enough,
\begin{equation*}
\prob\Big(\sup_D\phi_L\leq\Big(\sqrt{\frac{2}{\pi}}-\eta\Big)\ln\big(\sqrt{L}\big)\Big)\leq \exp\Big(-a\ln^2\big(\sqrt{L}\big)\Big).
\end{equation*}
\end{thm}

The proof of this theorem is essentially the same as the original. In the statement of Lemma \ref{sob}, one should introduce a compact $K\subset\mathring{\Sigma}$. The constant $C$ is then uniform for $p\in K$. In Proposition \ref{supup} the supremum is taken over $K$ and in the proof one should cover $K$ by small disks instead of $\Sigma$.

\section{Appendix}\label{sect_apx}

In section we recall some classical results from spectral theory used in the article and give an alternate characterization of the capacity which we use in Lemma \ref{huboff}.

\subsection{Classical spectral theory results}

In the following discussion we adopt the notations introduced in section \ref{sect_setting}. We begin by proving an approximation result for eigenfunctions of the laplacian.

\begin{lm}\label{bar_dens}
For any integers $h,k\geq 1$, the vector space spanned by the sequence $(\psi_n)_{n\geq 1}$ is dense in $C^h_0(\Sigma)$ and in $H^k_0(\Sigma)$.
\end{lm}
\begin{prf}
  By the Sobolev inequalities (see for instance, Theorem 5.6 (ii) of \cite{evans}), for any $h$ there exists $k_0$ such that for $k\geq k_0$, $H^k_0\subset C^h_0$ and the inclusion is continuous. On the other hand, if $h\geq k$, $C^h_0\subset H^k_0$ and the inclusion is continuous. Therefore, it suffices to prove that the space spanned by the $(\psi_n)_{n>0}$  is dense in all the $H^k_0$. By Garding's inequality (see Theorem (1), section 8 of \cite{taylor_ps}), for any $k\geq 1$, there is $C>0$ such that the scalar product
\[ (u\ |\ v)_k:= \langle u, (C+\Delta^k)v\rangle_2 \]
is equivalent to the standard $H^k_0$ scalar product. Let $f\in C^\infty_0$ and suppose that for any $n\geq 1$, $(f\ |\ \psi_n)_k=0$. Then
\[ 0=(f\ |\ \psi_n)_k=\langle f,(C+\Delta^k)\psi_n\rangle_2=(C+\lambda_n^k)\langle f,\psi_n\rangle_2. \]
Since $f\in L^2_0$ and $(\psi_n)_{n\geq 1}$ is dense in $L^2_0$, $f=0$. Since $C^\infty_0$ is dense in $H^k_0$, we conclude that so is $(\psi_n)_{n\geq 1}$.
\end{prf}

Now we discuss the functional calculus of the laplacian. This will be useful in the alternate characterization of the capacity. The operator $\Delta$ is a symmetric differential operator of order $2$ so it defines a bounded operator from $H^2(\Sigma)$ to $L^2(\Sigma)$ which, according to its spectral decomposition (see section \ref{sect_setting}) and by Lemma \ref{bar_dens}, defines an isomorphism $\Delta : H^2_0\rightarrow L^2_0$. Let $\Delta^{-1}$ be its inverse, which we extend to $L^2(\Sigma)$ by setting $\Delta^{-1} f=0$ for any constant function $f$. Since $\Delta$ is symmetric, $\Delta^{-1}$ is self-adjoint. Moreover, for any $n\geq 1$, $\Delta^{-1}\psi_n = \lambda_n^{-1}\psi_n$. For any $f\in L^2(\Sigma)$, we denote by $\sigma(f)$ the quadratic form $\sigma(f) := \langle f,\Delta^{-1}f\rangle_2$. By Parseval's formula,

\begin{equation}\eqlab \label{ref2}
\sigma(f) = \sum_{n\geq 1}\frac{1}{\lambda_n}\langle\psi_n,f\rangle_2^2.
\end{equation}
Finally, since $\Delta$ and $\Delta^{-1}$ are positive and $\Delta$ is a differential operator, then $\Delta:C^\infty_0\rightarrow C^\infty_0$ and $\Delta^{-1}:C^\infty_0\rightarrow C^\infty_0$ both admit square roots $\Delta^{1/2}$ and $\Delta^{-1/2}$ respectively, which are symmetric pseudo-differential operators of orders $1$ and $-1$ (see \cite{see_pow}). These define bounded operators
\begin{align*}
\Delta^{1/2} &: H^2\rightarrow H^1\\
\Delta^{1/2} &: H^1\rightarrow L^2\\
\Delta^{-1/2} &: L^2\rightarrow L^2.\\
\end{align*}
Finally, by construction, these operators restrict to $H^k_0$ and $L^2_0$ and for any $f\in H^2_0$, $\Delta^{-1/2}\Delta^{1/2}f=f$.

\subsection{The capacity}

For any subset $D\subset\Sigma$ that is not dense in $\Sigma$, there exist functions $h\in C^\infty_0(\Sigma)$ such that $h\geq 1$ on $D$. We define \textbf{the capacity} of $D$ the infimum over all such $h$ of $\frac{1}{2}\|\nabla h\|_2^2$ and denote it by $cap(D)$. This is a variation on the relative capacity used in \cite{bd_crit} and \cite{bdg_entrop}. Similarly to Lemma 2.1 of \cite{bd_crit}, in the case where $D$ is a non-empty proper open set with smooth boundary we have the following alternative characterization of the capacity.

\begin{prop}\label{capalt}
Let $D$ be a non-empty proper open subset of $\Sigma$ with smooth boundary. Then,
\[ cap(D)=\sup\Big\{\frac{1}{2\sigma(\mathbf{1}_Df)}\Big(\int_Df\Big)^2\ \Big|\ f\in C^\infty(\Sigma)\ f\geq0\text{ on $D$. }\Big\} \]
\end{prop}

In the proof of this proposition we will use the maximum principle and existence of solutions to certain elliptic PDEs. For this purpose, recall that we defined the laplacian to be $\Delta=d^*d$.

\begin{prf}[ of Proposition \ref{capalt}]
Let $h\in C_0^\infty(\Sigma)$ such that for each $x\in D$, $h(x)\geq 1$. Let $f\in C^\infty(\Sigma)$, non-negative on $D$.
\begin{align*}
\Big(\int_D f\Big)^2&\leq\Big(\int_D hf\Big)^2=\langle\mathbf{1}_Df,h\rangle_2^2=\langle\mathbf{1}_Df,\Delta^{-1/2}\Delta^{1/2}h\rangle_2^2=\langle\Delta^{-1/2}\mathbf{1}_Df,\Delta^{1/2}h\rangle_2^2\\
                    &\leq\langle\Delta^{-1/2}\mathbf{1}_Df,\Delta^{-1/2}\mathbf{1}_Df\rangle_2\langle\Delta^{1/2}h,\Delta^{1/2}h\rangle_2 \text{ by Cauchy-Schwarz}\\
                    &=\langle\mathbf{1}_Df,\Delta^{-1}\mathbf{1}_Df\rangle_2\langle h,\Delta h\rangle_2=\sigma(\mathbf{1}_Df)\|\nabla h\|^2_2.
\end{align*}
Passing to the infimum and supremum, since $D$ is non-empty, we deduce that
\[ \sup\Big\{\frac{1}{2\sigma(\mathbf{1}_Df)}\Big(\int_Df\Big)^2\ |\ f\in C^\infty(\Sigma)\ f\geq0\text{ on $D$. }\Big\}\leq cap(D). \]
We now show that one can find $f$ and $h$ so as to make the above functionals arbitrarily close to each other. The previous Cauchy-Schwarz inequality is close to equality when $\Delta^{-1}\mathbf{1}_Df$ is close to $h$. Recall that we defined $\Delta^{-1}$ to be zero when applied to constant functions. This means that we must find $f$ so that $\mathbf{1}_D f - \Delta h$ is close to a constant. We now define a would-be minimizer $h$. Let $\tau>0$ and $h$ be a continuous function constant equal to $1$ on $D$, smooth on $U=\Sigma\setminus\overline{D}$ such that $\Delta h$ is equal to $-\tau$ on $\Sigma\setminus\overline{D}$ (see Theorems 5 of section 6.2 and 3 of section 6.3 of \cite{evans} for the existence of $h$). Then $h$ is smooth up to the boundary of $U$ (see Theorem 6 of section 6.3 of \cite{evans}). We denote by $\partial_\nu h$ the outward normal of $h$, defined on $\partial U$ and by $d\sigma$ the volume density induced by $g$ on $\partial U=\partial D$. Since $\Delta h\leq 0$ on $U$, by the maximum principle, $h<1$ on $U$. Therefore, there exists $\tau$ such that $h$ has zero mean. Let $\Delta h$ be the laplacian of $h$ in the sense of distributions, which is well defined since $h$ is continuous. Before approximating $h$, we must find an explicit expression for $\Delta h$. For any $u\in C^\infty(\Sigma)$, by Stokes' theorem,
\begin{align*}
\langle \Delta h, u\rangle_2 = &\langle h, \Delta u\rangle_2 = \int_D h \Delta u |dV_g| + \int_U h \Delta u |dV_g|\\
                             = &\int_D \langle \nabla h, \nabla u\rangle_g |dV_g| + \int_{\partial D} h (-\partial_{\nu} u) d\sigma\\
                               &+ \int_U \langle \nabla h, \nabla u\rangle_g  |dV_g| + \int_{\partial U} h\partial_{\nu} u d\sigma\\
                             = &\int_U \langle \nabla h, \nabla u\rangle_g  |dV_g| = \int_U u \Delta h |dV_g| + \int_{\partial U}u\partial_\nu h d\sigma\\
                             = &\int_{\partial U}(\partial_\nu h) u d\sigma - \tau \int_U u |dV_g|.
\end{align*}
Therefore,
\[ \langle \Delta h +\tau, u\rangle_2 = \tau \int_D u |dV_g| + \int_{\partial(\Sigma\setminus D)}(\partial_\nu h) u d\sigma. \]
Again, by the maximum principle, $\partial_\nu h$ is positive. Let $(\tilde{h}_{\eps})_{\eps>0}$ be a sequence of smooth functions with zero mean converging uniformly to $h$ in $\Sigma$ that coincide with $h$ on $\Sigma\setminus D$. Then, for each $\eps>0$ there is $\delta(\eps)>0$ such that $h_\eps=(1+\delta(\eps))\tilde{h}_\eps\geq 1$ on $D$. Let $\tau_\eps = (1+\delta(\eps))\tau$. Then,  in the sense of distributions, $f_\eps=\Delta h_\eps + \tau_\eps$ is non-zero and non-negative on $\Sigma$, and vanishes on $\Sigma\setminus D$. Since it is smooth, it satisfies these properties in the classical sense as well. Moreover
\begin{align*}
(\int_Df_\eps)^2\geq&(\sup_{D}h_\eps)^{-2}(\int_D f_\eps h_\eps)^2\\
                   &=(\sup_{D}h_\eps)^{-2} \langle f_\eps, h_\eps \rangle_2^2=(\sup_{D}h_\eps)^{-2} \langle f_\eps, \Delta^{-1}f_\eps\rangle_{L^2}\langle h_\eps, \Delta h_\eps \rangle_2\\
                   &=(\sup_{D}h_\eps)^{-2}  \sigma(f_\eps)||\nabla h_\eps||_2^2=(\sup_{D}h_\eps)^{-2} \sigma(f_\eps\mathbf{1}_D)||\nabla h_\eps||_2^2.
\end{align*}
Therefore,
\[ cap(D)\leq \frac{1}{2}||\nabla h_\eps||_2^2\leq \frac{(\sup_{D}h_\eps)^2 }{2\sigma(f_\eps\mathbf{1}_D)}\Big(\int_Df_\eps\Big)^2. \]
Since $\sup_{D}h_\eps\xrightarrow[\eps\to 0]{} 1$, this concludes the proof of the lemma.
\end{prf}

\bibliography{mybib}{}
\bibliographystyle{plain}
\end{document}